\documentclass[final]{elsarticle_a}
\usepackage{amsmath}

\usepackage{amsfonts}

\usepackage{amssymb}

\usepackage{stmaryrd}

\usepackage[mathscr]{eucal}

\usepackage{amsbsy}

\usepackage{bm}

\usepackage{bbding}

\usepackage{array}

\usepackage[usenames]{color}

\usepackage{tabularx}

\usepackage{paralist}

\usepackage{fancyvrb}

\usepackage[position=below]{subfig}

\usepackage{boxedminipage}

\usepackage{multirow}

\usepackage{graphicx}
\usepackage{ifpdf}
\ifpdf
\DeclareGraphicsExtensions{.pdf,.png}
\else
\DeclareGraphicsExtensions{.eps}
\fi

\usepackage[draft,hypertexnames=false]{hyperref}

\usepackage[notref,notcite]{showkeys}

\newcommand{\Sec}[1]{\hyperref[sec:#1]{\S\ref*{sec:#1}}} 
\newcommand{\Section}[1]{\hyperref[sec:#1]{Section~\ref*{sec:#1}}} 
\newcommand{\App}[1]{\hyperref[sec:#1]{\ref*{sec:#1}}} 
\newcommand{\Eqn}[1]{\hyperref[eq:#1]{(\ref*{eq:#1})}} 
\newcommand{\Fig}[1]{\hyperref[fig:#1]{Figure~\ref*{fig:#1}}} 
\newcommand{\Tab}[1]{\hyperref[tab:#1]{Table~\ref*{tab:#1}}} 
\newcommand{\Thm}[1]{\hyperref[thm:#1]{Theorem~\ref*{thm:#1}}} 
\newcommand{\Cor}[1]{\hyperref[cor:#1]{Corollary~\ref*{cor:#1}}} 
\newcommand{\Alg}[1]{\hyperref[alg:#1]{Algorithm~\ref*{alg:#1}}} 
\newcommand{\Def}[1]{\hyperref[def:#1]{Definition~\ref*{def:#1}}} 

\newcommand{\Real}{{\mathbb R}}

\newenvironment{inlinemath}{$}{$}

\newcommand{\FMSFIG}[3]{
\begin{figure}[htbp]
  \centering
  \includegraphics[trim=0 30 0 15,clip]{#1}
  \caption{Comparison of accuracy for problems of size 
    $\mathbf{#2}$ for #3 with
    different amounts of missing data. The methods are CP-WOPT (red),
    INDAFAC (green), and EM-ALS (blue). We plot cumulative results for
    multiple starts. The first start is based on the $n$-mode 
    singular vectors, and the remaining starts are random.}
  \label{fig:#1}
\end{figure}
}

\newcommand{\TIMEFIG}[2]{
\begin{figure}[htbp]
  \centering
  \includegraphics[trim=0 25 0 0,clip]{#1}
  \caption{Comparison of timings for #2 for the total time for all starting points for
  different problem sizes and 
    different amounts of missing data. The methods are CP-WOPT-Dense
    (red),
    CP-WOPT-Sparse (magenta),
    INDAFAC (green), and EM-ALS (blue).}
  \label{fig:#1}
\end{figure}
}

\newcommand{\Tra}{^{{\sf T}}} 
\newcommand{\Tin}{^{-{\sf T}}} 

\newcommand{\V}[1]{{\bm{\mathbf{\MakeLowercase{#1}}}}} 
\newcommand{\Vbar}[1]{{\bm \bar{\mathbf{\MakeLowercase{#1}}}}} 
\newcommand{\VE}[2]{\MakeLowercase{#1}_{#2}} 
\newcommand{\VEP}[2]{\left( #1 \right)_{#2}} 
\newcommand{\Vn}[2]{\V{#1}^{(#2)}} 
\newcommand{\VnE}[3]{\MakeLowercase{#1}^{(#2)}_{#3}} 

\newcommand{\Oprod}{\circ} 

\newcommand{\M}[1]{{\bm{\mathbf{\MakeUppercase{#1}}}}} 
\newcommand{\Mhat}[1]{{\bm{\hat \mathbf{\MakeUppercase{#1}}}}} 
\newcommand{\Mbar}[1]{{\bm{\bar \mathbf{\MakeUppercase{#1}}}}} 
\newcommand{\ME}[2]{\MakeLowercase{#1}_{#2}} 
\newcommand{\MC}[2]{\V{#1}_{#2}} 
\newcommand{\Mn}[2]{\M{#1}^{(#2)}} 
\newcommand{\MnE}[3]{\MakeLowercase{#1}^{(#2)}_{#3}} 
\newcommand{\MnC}[3]{\V{#1}^{(#2)}_{#3}} 
\newcommand{\MbarnC}[3]{\Vbar{#1}^{(#2)}_{#3}} 
\newcommand{\MnCTra}[3]{\V{#1}^{(#2){{\sf T}}}_{#3}} 

\newcommand{\Khat}{\odot} 
\newcommand{\Hada}{\ast} 

\newcommand{\T}[1]{\boldsymbol{\mathscr{\MakeUppercase{#1}}}} 
\newcommand{\That}[1]{\boldsymbol{\hat \mathscr{\MakeUppercase{#1}}}} 
\newcommand{\Tbar}[1]{\boldsymbol{\bar \mathscr{\MakeUppercase{#1}}}} 
\newcommand{\TM}[2]{\M{#1}_{(#2)}} 
\newcommand{\TE}[2]{\MakeLowercase{#1}_{#2}} 
\newcommand{\TEP}[2]{\left( #1 \right)_{#2}} 

\newcommand{\Mz}[2]{\TM{#1}{#2}} 

\newcommand{\KOp}[1]{\llbracket #1 \rrbracket} 

\newcommand{\SizeN}[2]{{#1}_1 \times {#1}_2 \times \cdots \times {#1}_{#2}}

\newcommand{\SubscriptN}[2]{{#1}_1 {#1}_2 \cdots {#1}_{#2}}

\newcommand{\SumN}[3]{\sum_{{#1}_1=1}^{{#2}_1} %
  \sum_{{#1}_2=1}^{{#2}_2} %
  \cdots %
  \sum_{{#1}_{#3}=1}^{{#2}_{#3}}}

\newcommand{\norm}[1]{\left\lVert \, #1 \, \right\rVert}

\newcommand{\ip}[2]{\langle \, #1,#2 \, \rangle}

\newcommand{\FD}[2]{\frac{\partial #1}{\partial #2}}

\newcommand{\TheTitle}{Scalable Tensor Factorizations for Incomplete Data}

\newcommand{\TheAuthors}{Evrim Acar, Tamara G.\@ Kolda, Daniel M.\@ Dunlavy, and Morten M{\o}rup}

\newcommand{\TheKeywords}{missing data, incomplete data, tensor factorization, CANDECOMP, PARAFAC, optimization}

\begin{document}

  \title{\TheTitle\tnoteref{t1}}
  \author[tu]{Evrim Acar\fnref{fn1}}
  \ead{evrim.acar@bte.tubitak.gov.tr}

  \author[nm]{Daniel M. Dunlavy\fnref{fn1}}
  \ead{dmdunla@sandia.gov}

  \author[ca]{Tamara G. Kolda\corref{cor1}\fnref{fn1}}
  \ead{tgkolda@sandia.gov}

  \author[dk]{Morten M{\o}rup}
  \ead{mm@imm.dtu.dk}

  \cortext[cor1]{Corresponding author}
  \address[tu]{TUBITAK-UEKAE, Gebze, Turkey.}
  \address[ca]{Sandia National Laboratories, Livermore, CA 94551-9159.}
  \address[nm]{Sandia National Laboratories, Albuquerque, NM 87123-1318.}
  \address[dk]{Technical University of Denmark, 2800 Kgs.\@ Lyngby,
    Denmark.}
  \fntext[fn1]{This work was funded by the Laboratory Directed Research \&
    Development (LDRD) program at Sandia National Laboratories, a
    multiprogram laboratory operated by Sandia Corporation, a Lockheed
    Martin Company, for the United States Department of Energy's
    National Nuclear Security Administration under Contract
    DE-AC04-94AL85000.}
  \tnotetext[t1]{A preliminary conference version of this paper has
  appeared as \cite{SDM10}.}

\hypersetup{
 pdftitle={\TheTitle},
 pdfauthor={\TheAuthors},
 pdfkeywords={\TheKeywords}
}
 \hypersetup{colorlinks=true,urlcolor=black,citecolor=black,linkcolor=black}

\begin{abstract}
  The problem of incomplete data---i.e., data with missing or unknown
  values---in multi-way arrays is ubiquitous in
  biomedical signal processing, network traffic analysis,
  bibliometrics, social network analysis, chemometrics, computer
  vision, communication networks, etc.
  We consider the problem of how to factorize data sets with missing
  values with the goal of capturing the underlying
  latent structure of the data and possibly reconstructing missing
  values (i.e., tensor completion).
  We focus on one of the most well-known tensor factorizations that
  captures multi-linear structure, CANDECOMP/PARAFAC~(CP). In the
  presence of missing data, CP can be formulated as a weighted least
  squares problem that models \emph{only} the known entries. We
  develop an algorithm called CP-WOPT (CP Weighted OPTimization) that
  uses a first-order optimization approach to solve the weighted least
  squares problem.
  Based on extensive numerical experiments, our algorithm is shown to
  successfully factorize tensors with noise and up to 99\% missing
  data.  A unique aspect of our approach is that it scales to sparse
  large-scale data, e.g., $1000 \times 1000 \times 1000$ with five
  million known entries (0.5\% dense).
  We further demonstrate the usefulness of CP-WOPT on two real-world applications:
  a novel EEG (electroencephalogram) application
  where missing data is frequently encountered due to disconnections
  of electrodes and the problem of modeling computer network traffic
  where data may be absent due to the expense of the data collection process.
\end{abstract}


\maketitle

\section{Introduction}
\label{sec:introduction}

Missing data can arise in a variety of settings due to loss of information,
errors in the data collection process, or costly experiments.
For instance, in biomedical signal processing, missing data can be encountered during
EEG analysis, where multiple electrodes are used to collect the
electrical activity along the scalp. If one of the electrodes becomes loose
or disconnected, the signal is either lost or discarded due to
contamination with high amounts of mechanical noise. We also encounter the missing data
problem in other areas of data mining, such as packet losses
in network traffic analysis \cite{ZhRoWiQi09} and occlusions
in images in computer vision \cite{BuFi05}. Many real-world data with
missing entries are ignored because they are deemed unsuitable for
analysis, but this work contributes to the growing evidence that such
data can be analyzed.

\begin{figure}[tbph]
  \centering
  \includegraphics[width=1.5in,trim=10 10 100 10, clip]{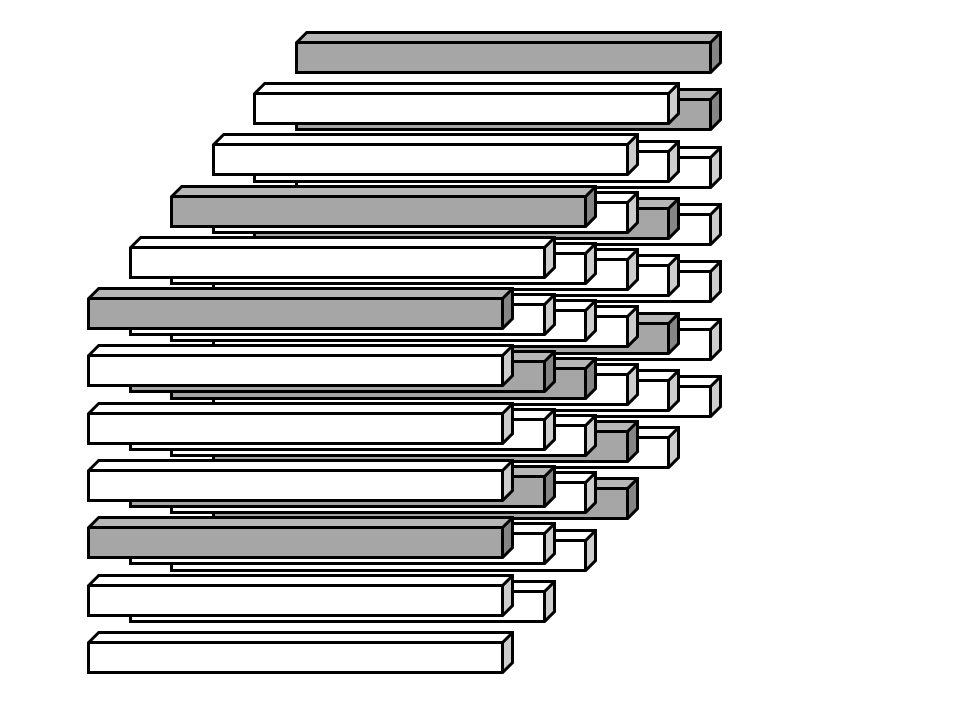}
  \caption[A $3$-way tensor with missing row fibers]%
  {A $3$-way tensor with missing row fibers (in gray).}
  \label{fig:fibers}
\end{figure}

Unlike most previous studies on missing data which have only considered matrices,
we focus here on the problem of missing data in \emph{tensors}
because it has been shown increasingly that data often have more
than two modes of variation and are therefore best represented as multi-way
arrays (i.e., tensors) \cite{AcYe09, KoBa09}. For instance, in EEG data
each signal from an electrode can be represented as a time-frequency matrix; thus,
data from multiple channels is three-dimensional (temporal, spectral, and
spatial) and forms a three-way array \cite{MiMaVaNi04}.
Social network data, network traffic data, and bibliometric data
are of interest to many applications such as community detection, link
mining, and more; these data can
have multiple dimensions/modalities, are often extremely large, and
generally have at least some missing data. 
These are just a few of the many data analysis applications where one needs to deal with large multi-way arrays with missing
entries. Other examples of multi-way arrays with missing entries from different disciplines
have also been studied in the literature \cite{ToBr05, Orekhov2003, Geng2009}. For instance, \cite{ToBr05}
shows that, in spectroscopy, intermittent machine failures or different sampling
frequencies may result in tensors with missing fibers (i.e.,
the higher-order analogues of matrix rows or columns, see \Fig{fibers}).
Similarly, missing fibers are encountered in multidimensional NMR
(Nuclear Magnetic Resonance) analysis, where
sparse sampling is used in order to
reduce the experimental time~\cite{Orekhov2003}.

Our goal is to capture the latent structure of the data via a
higher-order factorization, even in the presence of missing data.
Handling missing data in the context of matrix
factorizations, e.g., the widely-used principal component analysis,
has long been studied \cite{Ru74,GaZa79} (see \cite{BuFi05} for a review).
It is also closely related to the matrix completion problem, where the goal is to
recover the missing entries \cite{CaTa09,CaPl09}
(see \Sec{related} for more discussion). Higher-order factorizations, i.e., tensor factorizations,
have emerged as an important method for information analysis \cite{AcYe09, KoBa09}.
Instead of flattening (unfolding) multi-way arrays into matrices
and using matrix factorization techniques, tensor models
preserve the multi-way nature of the data and extract
the underlying factors in each mode (dimension) of a higher-order array.

We focus here on the CANDECOMP/PARAFAC (CP) tensor decomposition
\cite{CaCh70,Ha70}, which is a tensor model commonly used in various
applications \cite{MiMaVaNi04,KoBaKe05,Br06a,AcBiBiBr07,MoHaAr07}.
To illustrate differences between matrix and tensor factorizations, we introduce the CP decomposition for three-way tensors; discussion of the CP decomposition for general $N$-way tensors can be found in \Sec{algorithm}. Let $\T{X}$ be a three-way tensor of size $I
\times J \times K$, and assume its rank is $R$
(see \cite{KoBa09} for a detailed discussion on tensor rank).
With perfect data,
the CP decomposition is defined by \emph{factor matrices} $\M{A}$,
$\M{B}$, and $\M{C}$ of sizes $I \times R$, $J \times R$, and $K
\times R$, respectively, such that
\begin{displaymath}
  \TE{x}{ijk} = \sum_{r=1}^R \ME{A}{ir} \ME{B}{jr} \ME{C}{kr},
  \quad\text{for all } i = 1,\dots,I, \; j = 1,\dots,J, \; k = 1,\dots,K.
\end{displaymath}

\begin{figure}[t]
\centering
\includegraphics[width=3.3in]{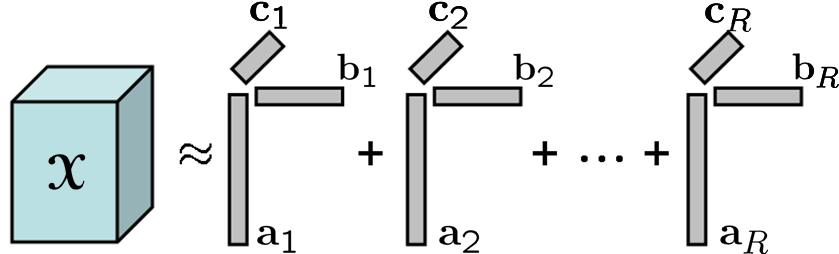}
\caption[Illustration of an $R$-component CP model for a
third-order tensor]{Illustration of an $R$-component CP model for a
third-order tensor $\T{X}$.} \label{fig:CP}
\end{figure}

In the presence of noise, the true $\T{X}$ is not observable and we
cannot expect equality. Instead, the CP decomposition should minimize
the error function
\begin{equation}
  \label{eq:fABC}
  f(\M{A},\M{B},\M{C}) =
  \frac{1}{2}
  \sum_{i=1}^I \sum_{j=1}^J \sum_{k=1}^K \left( \TE{x}{ijk} -
  \sum_{r=1}^R \ME{A}{ir} \ME{B}{jr} \ME{C}{kr} \right)^2 .
\end{equation}
An illustration of CP for third-order tensors is given in \Fig{CP}.
The CP decomposition is extensible to $N$-way tensors for $N \geq
3$, and there are numerous methods for computing it \cite{AcKoDu09}.

In the case of incomplete data, a standard practice is to impute the missing
values in some fashion (e.g., replacing the missing entries using average
values along a particular mode). Imputation can be useful as long as the amount of missing
data is small; however, performance degrades for large amounts of
missing data \cite{Ru74,SDM10}. As a better alternative,
factorizations of the data with imputed values for missing entries can be used to 
re-impute the missing values
and the procedure can be repeated to iteratively determine suitable values
for the missing entries. Such a procedure is an example of the
expectation maximization (EM) algorithm \cite{DeLaRu77}.
Computing CP decompositions by combining the alternating least squares method,
which computes the factor matrices one at a time, and iterative imputation
(denoted EM-ALS in this paper) has been shown to be
quite effective and has the advantage of often being simple and fast.
Nevertheless, as the amount of missing data increases, the performance
of the algorithm may suffer since the initialization and the
intermediate models used to impute the missing values will increase
the risk of converging to a less than optimal factorization \cite{ToBr05}. Also, the
poor convergence of alternating methods due to their vulnerability to
flatlining, i.e., stagnation, is noted in \cite{BuFi05}.

In this paper, though, we focus on
using a weighted version of the error function
to ignore missing data and model only the known entries.
In that case, nonlinear optimization can be used to directly solve the
weighted least squares
problem for the CP model. The weighted version of \Eqn{fABC} is

\begin{equation}
  \label{eq:fWABC}
  f_{\T{W}}(\M{A},\M{B},\M{C}) =
  \frac{1}{2}
  \sum_{i=1}^I \sum_{j=1}^J
  \sum_{k=1}^K \left\{ \TE{w}{ijk} \left( \TE{x}{ijk} -
  \sum_{r=1}^R \ME{A}{ir} \ME{B}{jr} \ME{C}{kr} \right) \right\}^2 ,
\end{equation}
where $\T{W}$, which is the same size as $\T{X}$, is a nonnegative weight tensor defined as
\begin{displaymath}
  \TE{W}{ijk} =
  \begin{cases}
    1 & \text{if $\TE{X}{ijk}$ is known}, \\
    0 & \text{if $\TE{X}{ijk}$ is missing}, \\
  \end{cases}
  \quad\text{for all } i = 1,\dots,I, \; j = 1,\dots,J, \; k = 1,\dots,K.
\end{displaymath}

Our contributions in this paper are summarized as follows.
\begin{inparaenum}[(a)]
\item We develop a scalable algorithm called CP-WOPT (CP Weighted OPTimization)
    for tensor factorizations in the presence of missing data. CP-WOPT
    uses first-order optimization to solve the weighted least squares
    objective function over all the factor matrices simultaneously.
\item We show that CP-WOPT can scale to sparse, large-scale data
  using specialized sparse data structures, significantly reducing the storage and computation
  costs.
\item Using extensive numerical experiments on simulated data sets, we
  show that CP-WOPT can successfully factor tensors with noise and up
  to 99\% missing data. In many cases, CP-WOPT is significantly faster
  than the best published direct optimization method in the literature \cite{ToBr05}.
\item We demonstrate the applicability of the proposed algorithm on a
  real data set in a novel EEG application where data is incomplete
  due to failures of particular electrodes. This is a common occurrence
  in practice, and our experiments show that even if signals from
  almost half of the channels are missing, underlying brain activities
  can still be captured using the CP-WOPT algorithm, illustrating the
  usefulness of our proposed method. 
\item In addition to tensor factorizations, we also show that CP-WOPT can
  be used to address the tensor completion problem in the context of network
  traffic analysis. We use the factors captured by the CP-WOPT algorithm to reconstruct
  the tensor and illustrate that even if there is a large amount of missing data, the
  algorithm is able to keep the relative error in the missing entries close to the modeling error.
\end{inparaenum}

The paper is organized as follows. We introduce the notation used throughout
the paper in \Sec{notation}.
In \Sec{related}, we discuss related work in matrix and tensor factorizations.
The computation of the function and gradient values for the
general $N$-way weighted version of the error function
and the presentation of the CP-WOPT method are given in
\Sec{algorithm}. Numerical results on both simulated and real data are
given in \Sec{experiments}. Conclusions and future work are discussed
in \Sec{conclusions}.

\section{Notation}
\label{sec:notation}

Tensors of order $N \geq 3$ are denoted by Euler script letters
($\T{X},\T{Y},\T{Z}$),
matrices are denoted by boldface capital letters
($\M{A},\M{B},\M{C}$),
vectors are denoted by boldface lowercase letters
($\V{A},\V{B},\V{C}$),
and scalars are denoted by lowercase letters ($a$, $b$, $c$).
Columns of a matrix are denoted by boldface lower letters with a
subscript ($\MC{A}{1},\MC{A}{2},\MC{A}{3}$ are first three columns of
$\M{A}$).
Entries of a matrix or a tensor are denoted by lowercase letters with
subscripts, i.e., the $(i_1, i_2, \dots, i_N)$ entry of an $N$-way
tensor $\T{X}$ is denoted by $\TE{X}{\SubscriptN{i}{N}}$.

An $N$-way tensor can be rearranged as a matrix; this is called \emph{matricization},
also known as \emph{unfolding} or \emph{flattening}. The mode-$n$ matricization of a tensor
$\T{X} \in \Real^{\SizeN{I}{N}}$ is denoted by $\Mz{X}{n}$ and
arranges the mode-$n$ one-dimensional ``fibers'' to be the columns of
the resulting matrix; see \cite{SDM10,KoBa09} for details.

Given two tensors $\T{X}$ and $\T{Y}$ of equal size $\SizeN{I}{N}$,
their Hadamard (elementwise) product is denoted by $\T{X} \Hada \T{Y}$
and defined as
\begin{displaymath}
  \TEP{\T{X} \Hada \T{Y}}{\SubscriptN{i}{N}} =
  \TE{X}{\SubscriptN{i}{N}} \TE{Y}{\SubscriptN{i}{N}}
  \quad\text{ for all } i_n \in \{1,\dots,I_n\} \text{ and } n \in \{1,\dots,N\}
\end{displaymath}
The \emph{inner product} of two same-sized tensors $\T{X},\T{Y} \in
\Real^{\SizeN{I}{N}}$ is the sum of the products of their entries, i.e.,
\begin{displaymath}
\ip{\T{X}}{\T{Y}}
=
\SumN{i}{I}{N}
\TE{X}{\SubscriptN{i}{N}}
\TE{Y}{\SubscriptN{i}{N}}.
\end{displaymath}
For a tensor $\T{X}$ of size $\SizeN{I}{N}$, its \emph{norm} is
\begin{inlinemath}
  \norm{ \T{X} } = \sqrt{\ip{\T{X}}{\T{X}}}.
\end{inlinemath}
For matrices and vectors, $\|\cdot\|$ refers to the analogous
Frobenius and two-norm, respectively.
We can also define a weighted norm as follows. Let $\T{X}$ and
$\T{W}$ be two tensors of size $\SizeN{I}{N}$. Then the
$\T{W}$-weighted norm of $\T{X}$ is
\begin{displaymath}
  \norm{ \T{X} }_{\T{W}} = \norm{ \T{W} \Hada \T{X} }.
\end{displaymath}

Given a sequence of matrices $\Mn{A}{n}$ of size $I_n \times R$ for
$n=1,\dots,N$, the notation
$\KOp{\Mn{A}{1},\Mn{A}{2},\dots,\Mn{A}{N}}$ defines an $\SizeN{I}{N}$
tensor whose elements are given by
\begin{displaymath}
  \TEP{\KOp{\Mn{A}{1},\Mn{A}{2},\dots,\Mn{A}{N}}}{\SubscriptN{i}{N}}
  = \sum_{r=1}^R \prod_{n=1}^N \MnE{a}{n}{i_nr},
 \text{ for }
 i_n \in \{1,\dots,I_n\}, n \in \{1,\dots,N\}.
\end{displaymath}
For just two matrices, this reduces to familiar expressions:
\begin{inlinemath}
  \KOp{\M{A},\M{B}} = \M{A}\M{B}\Tra.
\end{inlinemath}
Using the notation defined here, \Eqn{fWABC} can be rewritten as
\begin{displaymath}
    f_{\T{W}}(\M{A},\M{B},\M{C}) =
    \frac{1}{2}
    \norm{
      \T{X} - \KOp{\M{A},\M{B},\M{C}}
    }_{\T{W}}^2.
\end{displaymath}

\section{Related Work in Factorizations with Missing Data}
\label{sec:related}

In this section, we first review the approaches for handling missing
data in matrix factorizations and then discuss how these techniques
have been extended to tensor factorizations.

\subsection{Matrix Factorizations}

Matrix factorization in the presence of missing entries is a problem
that has been studied for several decades; see, e.g.,
\cite{Ru74,GaZa79}. The problem is typically formulated analogously to
\Eqn{fWABC} as
\begin{equation}
  \label{eq:fWAB}
  f_{\M{W}} (\M{A},\M{B})
  = \frac{1}{2} \norm{
  \M{X} - \M{A}\M{B}\Tra
  }_{\M{W}}^2.
\end{equation}
A common procedure for solving this problem is EM which combines
imputation and alternation 
 \cite{ToBr05,SrJa03}. In this approach, the missing values of $\M{X}$ are imputed using the
current model, $\Mhat{X} = \M{A}\M{B}\Tra$, as follows:
\begin{displaymath}
  \Mbar{X} = \M{W} \Hada \M{X} + (\M{1} - \M{W}) \Hada \Mhat{X},
\end{displaymath}
where $\M{1}$ is the matrix of all ones.
Once $\Mbar{X}$ is generated, the matrices $\M{A}$ and $\M{B}$ can
then be alternatingly updated according to
the error function
\begin{inlinemath}
 \frac12 \| \Mbar{X} - \M{A}\M{B}\Tra \|^2
\end{inlinemath}
(e.g., using the linear least squares method).
See \cite{SrJa03, Ki97} for further discussion of the EM method in the missing data and
general weighted case.

Recently, a direct nonlinear optimization approach was proposed for
matrix factorization with missing data~\cite{BuFi05}. In this case,
\Eqn{fWAB} is solved directly using a 2nd-order damped Newton
method. This new method is compared to other standard techniques based
on some form of alternation and/or imputation as well as hybrid
techniques that combine both approaches. Overall, the conclusion is
that nonlinear optimization strategies are key to successful matrix
factorization. Moreover,  the authors observe that the alternating methods tend to
take much longer to converge to the solution even though they make
faster progress initially.
This work is theoretically the most related to what we propose---the main
differences are 1) we focus on tensors rather than matrices, and
2) we use first-order rather than second-order optimization methods
(we note that first order methods are mentioned as future work in \cite{BuFi05}).

A major difference between matrix and tensor factorizations is worth
noting here. In \cite{SrJa03,BuFi05}, the lack of uniqueness in matrix
decompositions is discussed. Given any invertible matrix $\M{G}$,
\begin{inlinemath}
  \KOp{\M{A},\M{B}} = \KOp{\M{AG}, \M{BG}\Tin}.
\end{inlinemath}
This means that there is an infinite family of equivalent solutions to \Eqn{fWAB}.
In \cite{BuFi05}, regularization is recommended as a partial solution;
however regularization can only control scaling indeterminacies and not rotational freedom. In the
case of the CP model, often a unique solution (including
trivial indeterminacies of scaling and column permutation) can be
recovered exactly; see, e.g., \cite{KoBa09} for further
discussion on uniqueness of the CP decomposition.

Factorization of matrices with missing entries is also closely related to the matrix
completion problem. In matrix completion, one tries to recover the missing
matrix entries using the low-rank structure of the matrix. Recent work in
this area \cite{CaTa09,CaPl09} shows that even if a small amount of matrix entries
are available and those are corrupted with noise, it is still possible to recover the
missing entries up to the level of noise. In \cite{CaPl09}, it is also discussed how this problem relates
to the field of compressive sensing, which exploits
structures in data to generate more compact representations of the data. Practically speaking, the difference between completion and factorization
is how success is measured. Factorization methods aim to increase accuracy in
the factors; in other words, capture the underlying phenomena as well as possible. 
Completion methods, on the other hand, seek accuracy in filling in the missing data.
Obviously, once a factorization has been computed, it can be used to reconstruct
the missing entries. In fact, many completion methods use this procedure.

\subsection{Tensor Factorizations}

The EM procedure discussed for matrices has also been
widely employed for tensor factorizations with missing data.
If the current model is
$\KOp{\M{A},\M{B},\M{C}}$, then we fill in the missing entries of
$\T{X}$ to produce a complete tensor according to
\begin{displaymath}
  \Tbar{X} = \T{W} \Hada \T{X}
  + (\M{1} - \T{W}) \Hada \KOp{\M{A},\M{B},\M{C}},
\end{displaymath}
where $\M{1}$ is the tensor of all ones the same size as $\T{W}$. The factor matrices are then
updated using alternating least squares (ALS) as those that best fit $\Tbar{X}$.
See, e.g., \cite{Br98,WaMa01} for further details. As noted
previously, we denote this method as EM-ALS.

Paatero \cite{Pa97} and Tomasi and Bro \cite{ToBr05} have
investigated direct nonlinear approaches based on Gauss-Newton
(GN). The code from \cite{Pa97} is not widely available; therefore,
we focus on \cite{ToBr05} and its INDAFAC (INcomplete DAta
paraFAC) procedure which specifically uses the Levenberg-Marquardt
version of GN for fitting the CP model to data with missing entries.
The primary application in \cite{ToBr05} is missing
data in chemometrics experiments.  This approach is compared to EM-ALS
with the result being that INDAFAC and EM-ALS perform almost equally
well in general with the exception that INDAFAC is more accurate
for difficult problems, i.e., higher collinearity and systematically
missing patterns of data. In terms of computational
efficiency, EM-ALS is usually faster but becomes slower
than INDAFAC as the percentage of missing entries increases and also
depending on the missing entry patterns.

Both INDAFAC and CP-WOPT address the problem of fitting the CP
model to incomplete data sets by solving \Eqn{fWABC}. The difference
is that
INDAFAC is based on second-order optimization while CP-WOPT is first-order
with a goal of scaling to larger problem sizes.
\section{CP-WOPT Algorithm}
\label{sec:algorithm}

We consider the general $N$-way CP factorization problem for tensors
with missing entries. Let $\T{X}$ be a real-valued tensor of size $\SizeN{I}{N}$
and assume its rank is known to be $R$.\footnote{In practice, the rank is
  generally not known and is not easily determined. Understanding the
  performance of the methods under consideration in that scenario is a
  topic of future work.
  Results in \cite{AcKoDu09} indicate that direct optimization methods
  have an advantage over alternating least squares approaches
  when the rank is overestimated.}  Define a
nonnegative weight tensor $\T{W}$ of the same size as $\T{X}$ such
that
\begin{multline}
  \label{eq:W}
  \TE{W}{\SubscriptN{i}{N}} =
  \begin{cases}
    1 & \text{if $\TE{X}{\SubscriptN{i}{N}}$ is known},\\
    0 & \text{if $\TE{X}{\SubscriptN{i}{N}}$ is missing},
  \end{cases}\\
  \quad \text{for all } i_n \in \{1,\dots,I_n\} \text{ and } n \in \{1,\dots,N\}.
\end{multline}
The $N$-way objective function is defined by
\begin{equation}
\label{eq:fWNway}
   f_{\T{W}}(\Mn{A}{1},\Mn{A}{2},\dots,\Mn{A}{N})
      =  \frac{1}{2}
      \norm{ \left( \T{X} - \KOp{\Mn{A}{1},\dots,\Mn{A}{N}} \right)}_{\T{W}}^2 .
\end{equation}
It is a function of $P =R \sum_{n=1}^N I_n$ variables.
Due to the well-known indeterminacies of the CP model, it may also be
desirable to add regularization to the objective function as in
\cite{AcKoDu09}, but this has not been necessary thus far in our experiments.

Equation \Eqn{fWNway} is equivalent to
\begin{multline}
  \label{eq:YZ}
     f_{\T{W}}(\Mn{A}{1},\Mn{A}{2},\dots,\Mn{A}{N}) =
     \frac{1}{2} \norm{\T{Y} - \T{Z}}^2, \\
     \quad\text{where}\quad
     \T{Y} = \T{W} \Hada \T{X} \text{ and }
     \T{Z} = \T{W} \Hada \KOp{\Mn{A}{1},\dots,\Mn{A}{N}}.
\end{multline}
The tensor $\T{Y}$ can be precomputed as neither $\T{W}$ nor $\T{X}$
change during the iterations. We will see later that $\T{Z}$ can be
computed efficiently for sparse $\T{W}$.

Our goal is to find matrices $\Mn{A}{n} \in \Real^{I_n \times R}$ for
$n=1,\dots,N$ that minimize the weighted objective function in
\Eqn{fWNway}.
We 
show how to compute the function and gradient values efficiently when
$\T{W}$ is dense in \Sec{dense} and when
$\T{W}$ is sparse in \Sec{sparse}. Once the function and gradient are
known, any gradient-based optimization method \cite{NoWr99} can be
used to solve the optimization problem.

The derivation of the gradient in the weighted case is given in
\cite{SDM10}; here we just report the formula.
In matrix notation, using $\T{Y}$ and $\T{Z}$ from \Eqn{YZ}, we can
rewrite the gradient equation as
\begin{align}
  \label{eq:grad}
  \FD{f_{\T{W}}}{\Mn{A}{n}} = \left( \Mz{Z}{n} - \Mz{Y}{n} \right) \Mn{A}{-n},
\end{align}
where
\begin{displaymath}
\Mn{A}{-n} = \Mn{A}{N} \Khat \cdots \Khat \Mn{A}{n+1} \Khat
\Mn{A}{n-1} \Khat \cdots \Khat \Mn{A}{1}
\end{displaymath}
for $n=1,\dots,N.$ The symbol $\Khat$ denotes the Khatri-Rao product;
see, e.g., \cite{SDM10,KoBa09} for details.

\subsection{Computations with Dense $\T{W}$}
\label{sec:dense}

\Fig{code} shows the algorithmic steps for computing the function and
gradient values. 
We are optimizing a function of $P$
variables as defined by the factor matrices $\Mn{A}{1}$ through
$\Mn{A}{N}$. The gradient is computed as a series of
matrices $\Mn{G}{n} \equiv \FD{f}{\Mn{A}{n}}$ for $n=1,\dots,N$.
Recall that $\Mz{T}{n}$ in the last step is the
unfolding of the tensor $\T{T}$ in mode $n$.


\subsection{Computations with Sparse $\T{W}$}
\label{sec:sparse}

If a large number of entries is missing, then $\T{W}$ is
sparse. In this case, there is no need to allocate storage for every
entry of the
tensor $\T{X}$. Instead, we can  store and work with just the known
values, making the method efficient in both storage and time.
\Fig{sparse} shows analogous computations to \Fig{code} in the case
that $\T{W}$ is sparse.

\begin{figure}[tbp]
  \centering
  \subfloat[Dense $\T{W}$.]{
  \label{fig:code}
  \begin{boxedminipage}{\textwidth}
    \begin{enumerate}
    \item Assume $\T{Y} = \T{W} \Hada \T{X}$ and $\gamma =
      \norm{\T{Y}}^2$ are precomputed.
    \item Compute $\T{z} = \T{W} \Hada
      \KOp{\Mn{A}{1},\dots,\Mn{A}{N}}$.
    \item Compute function value: $f = \frac{1}{2} \gamma - \ip{ \T{Y}}{\T{Z}
      } + \frac{1}{2} \norm{\T{Z}}^2$.
    \item Compute $\T{T} = \T{Y} - \T{Z}.$
    \item Compute gradient matrices: $\Mn{G}{n} = -\Mz{T}{n} \Mn{A}{-n}$ for $n =
      1,\dots,N$.
    \end{enumerate}
  \end{boxedminipage}
  }\\
  \subfloat[Sparse $\T{W}$.]{
  \label{fig:sparse}
  \begin{boxedminipage}{\textwidth}
    \begin{enumerate}
    \item Let $\mathcal{I} = \{\Vn{i}{1},\Vn{i}{2},\dots,\Vn{i}{Q}\}$ be
      an ordered set of all the locations where $\T{W}=1$, i.e., the
      indices of the known values.
      Let $\V{y}$ be the length-$Q$ vector of the values of
      $\T{X}$ at the locations indicated by $\mathcal{I}$. Assume
      $\V{y}$ and $\gamma = \norm{\V{y}}^2$ are precomputed.
  \item Compute the $Q$-vector $\V{z}$ as
    \begin{displaymath}
      \VE{z}{q} = \sum_{r=1}^R \prod_{n=1}^N
      \MnE{A}{n}{\VnE{i}{q}{n} r}
      \quad\text{for}\quad q = 1,\dots,Q.
    \end{displaymath}
  \item Compute function value: $f = \frac{1}{2} \gamma - \V{y}\Tra\V{z} +
    \frac{1}{2} \norm{\V{z}}^2$.
  \item Compute $\V{t} = \V{y} - \V{z}$.
  \item Compute gradient matrices $\Mn{G}{n}$ for $n=1,\dots,N$ as follows:
    \begin{displaymath}
      \MnE{G}{n}{jr} = -\sum_{{q=1}\atop{q:\VnE{i}{q}{n} = j}}^Q
      \left( \VE{t}{q}
      \prod_{{m=1}\atop{m \neq n}}^N \MnE{a}{m}{\VnE{i}{q}{m} r}
      \right).
    \end{displaymath}
  \end{enumerate}
  \end{boxedminipage}}
  \caption{CP-WOPT computation of function value and gradient.}
\end{figure}

Let the ordered set
$\mathcal{I} = \{\Vn{i}{1},\Vn{i}{2},\dots,\Vn{i}{Q}\}$ be the indices of the known
values, i.e., all the locations where $\T{W}=1$. Each $\Vn{i}{q}$, $q \in \{1, \dots, Q\}$, is
an N-tuple of indices whose $n$th entry is denoted by $\VnE{i}{q}{n}$.
The known entries of $\T{X}$ can be stored in an array $\V{y}$ of length $Q$
so that
\begin{displaymath}
  \VE{y}{q} = \TE{X}{\VnE{i}{q}{1}, \VnE{i}{q}{2}, \dots,
    \VnE{i}{q}{N}}
  \quad\text{for}\quad q=1,\dots,Q.
\end{displaymath}
Observe that
\begin{inlinemath}
 \norm{\T{W}\Hada\T{X}}^2 = \norm{\V{y}}^2.
\end{inlinemath}
Thus, the value of $\gamma$ is the same in both \Fig{code} and \Fig{sparse}.

In the dense version (\Fig{code}), we need to compute $\T{z} =
\T{W} \Hada \KOp{\Mn{A}{1},\dots,\Mn{A}{N}}$. In the sparse version,
we observe that $\T{z}$ is necessarily zero anywhere there is a
missing entry in $\T{X}$. Consequently, it is only necessary to
compute the entries of
$\T{z}$ corresponding to known values in $\mathcal{I}$; the vector $\V{z}$
corresponds to these known entries.
The computations for Step 2 in
\Fig{sparse} can be done efficiently in MATLAB using the
``expanded'' vectors technique of \cite[\S3.2.4]{BaKo07}.
Let the $r$th summand be denoted by $\V{u}$, i.e.,
\begin{displaymath}
  \VE{u}{q} = \prod_{n=1}^N
      \MnE{A}{n}{\VnE{i}{q}{n} r}
      \quad\text{for}\quad q = 1,\dots,Q.
\end{displaymath}
Let $\Vn{v}{n}$ be ``expanded'' vectors
of length $Q$ for $n=1,\dots,N$ defined by
\begin{equation}
  \label{eq:vn}
  \VnE{v}{n}{q} = \MnE{A}{n}{\VnE{i}{n}{q} r}
  \quad\text{for}\quad
  q =1,\dots,Q.
\end{equation}
Then the vector $\V{u}$ can be calculated as
\begin{displaymath}
  \V{u} = \Vn{V}{1} \Hada \Vn{V}{2} \Hada \cdots \Hada \Vn{V}{N}.
\end{displaymath}
This can naturally be done iteratively (i.e., only one $\Vn{v}{n}$ is
calculated at a time) to reduce storage costs. Likewise, each summand
$\V{u}$ can be iteratively added to compute the final answer $\V{z}$.

In Step 3, the function values in \Fig{code} and \Fig{sparse} are clearly
equivalent since $\V{Y}$ and $\V{z}$ contain the
nonzero values of $\T{Y}$ and $\T{z}$, respectively.  Similarly, the
vector $\V{t}$ in Step 4 of \Fig{sparse} represents just the nonzero
values of $\T{T}$ in \Fig{code}.

The computation of the gradients in Step 6 of
\Fig{sparse} performs a matricized-tensor-times-Khatri-Rao-product
(mttkrp) calculation, which has been described for the sparse data
case in \cite[\S5.2.6]{BaKo07}. Here we briefly summarize the
methodology. The $r$th column of $\Mn{G}{n}$, $\MnC{G}{n}{r}$ is calculated as
follows. Let the vectors $\Vn{v}{n}$ be defined as above in
\Eqn{vn}, but define $\V{u}$ instead as
\begin{displaymath}
  \V{u} = \V{t} \Hada \Vn{V}{1} \Hada \Vn{V}{2} \Hada \cdots \Vn{v}{n-1}
  \Hada \Vn{v}{n+1} \Hada \cdots \Hada \Vn{V}{N}.
\end{displaymath}
Then
\begin{displaymath}
  \VEP{\MnC{G}{n}{r}}{j} = \sum_{\VnE{i}{q}{n}=j}
  \VE{u}{q}
  \quad\text{for}\quad j =1,\dots,I_n.
\end{displaymath}
This can be computed efficiently using the \texttt{accumarray}
function in MATLAB, and the code is available in version 2.4 of the
Tensor Toolbox for MATLAB \cite{TTB}.


\section{Experiments}
\label{sec:experiments}

On both real and simulated three-way data, we assess the performance
of the CP-WOPT method. We compare CP-WOPT with other methods and also
demonstrate its performance on two applications.

\subsection{Computational environment}
\label{sec:env}
All experiments were performed using MATLAB 2009b on a Linux
Workstation (RedHat 5.2) with 2 Quad-Core Intel Xeon 3.0GHz processors
and 32GB RAM. Timings were performed using MATLAB's \texttt{tic} and
\texttt{toc} functions since \texttt{cputime} is known to produce
inaccurate results for multi-CPU and/or multi-core systems.

CP-WOPT is implemented in the Tensor Toolbox \cite{TTB}. We
consider dense and sparse versions, based on the gradient and function
computations shown in \Fig{code} and \Fig{sparse}, respectively. We
use the nonlinear conjugate gradient (NCG) method with Hestenes-Stiefel updates \cite{NoWr99} and the
Mor\'{e}-Thuente line search \cite{MoTh94} provided in the Poblano Toolbox
\cite{DuKoAc10} as the optimization method.

We compare CP-WOPT to two other methods: EM-ALS (implemented in the N-way
Toolbox for MATLAB, version 3.10 \cite{AnBr00}) and INDAFAC \cite{INDAFAC},
which is a damped Gauss-Newton method proposed by Tomasi and Bro \cite{ToBr05}.
Previously, Tomasi and Bro showed that INDAFAC
converged to solutions in many fewer iterations than EM-ALS.

The stopping conditions are set as follows.
All algorithms use the relative change in the function value
$f_{\T{W}}$ in \Eqn{fWABC} as a stopping condition (set to $10^{-8}$).
In INDAFAC, the tolerance on the infinity norm of
the gradient is set to $10^{-8}$ and the maximum number of iterations
is set to $500$.
In CP-WOPT, the tolerance on the two-norm of the gradient divided by the
number of entries in the gradient is set to $10^{-8}$, the maximum
number of iterations is set to $500$, and the maximum number of
function evaluations is set to $10000$.  In EM-ALS, the maximum number
of iterations (equivalent to one function evaluation) is set to 10000.

All the methods under consideration are iterative methods.  We used
multiple starting points for each randomly generated problem.
The first starting point is generated using the $n$-mode singular
vectors \footnote{The $n$-mode singular vectors are the left
singular vectors of $\Mz{X}{n}$ ($\T{X}$ unfolded in mode $n$).} of $\T{X}$ with
missing entries replaced by zero. In our preliminary experiments, this
starting procedure produced significantly better results than random
initializations, even with large amounts of missing data. In order to
improve the chances of reaching the global minimum, 
additional starting points were generated randomly. The same set of
starting points was used by all the methods in the same order.

\subsection{Validation metrics}

If the true factors are known, then we can assess the recovery of the
factors via the \emph{factor match score} (FMS) defined as follows.
Let the correct and computed factorizations be given by
\begin{displaymath}
  \sum_{r=1}^R \lambda_r \; \MnC{A}{1}{r} \Oprod \MnC{A}{2}{r}
  \Oprod \cdots \Oprod \MnC{A}{N}{r}
  \quad\text{and}\quad
  \sum_{r=1}^{\bar R} \bar\lambda_r \; \MbarnC{A}{1}{r} \Oprod \MbarnC{A}{2}{r}
  \Oprod \cdots \Oprod \MbarnC{A}{N}{r},
\end{displaymath}
respectively. Without loss of generality, we assume that all the
vectors have been scaled to unit length and that the
scalars are positive. Further, we assume $\bar R
\geq R$ so that the computed solution has at least as many components
as the true solution. (One could add all-zero components to the
computed solution if $\bar R < R$.) Recall that there is a permutation
ambiguity, so all possible matchings of
components between the two solutions must be considered.
Under these conditions, the FMS is defined as
\begin{equation}
\label{eq:FMS}
  \text{FMS} = \max_{\sigma\in\Pi(R,\bar R)}
  \frac{1}{R} \sum_{r=1}^R
  \left(
    1 - \frac{|\lambda_r -  \bar\lambda_{\sigma(r)}|}
    {\max\{\lambda_r,\bar\lambda_{\sigma(r)}\}}
  \right)
  \prod_{n=1}^N |\MnCTra{A}{n}{r} \MbarnC{A}{n}{\sigma(r)}| .
\end{equation}
If $\bar R = R$, then the set $\Pi(R,\bar R)$ is all
permutations of 1 to $R$; otherwise, it is all possible permutations of all
$\bar R \choose R$ mappings of $\{1,\dots,\bar R\}$ to
$\{1,\dots,R\}$. The FMS can be between 0 and 1, and the best possible
FMS is 1.
If $\bar R \geq R$, some components in the computed solution are
completely ignored in the calculation of FMS, but this is not an issue
for us because we use $\bar R = R$ throughout.

We also consider the problem of recovering missing data. Let $\T{X}$
be the original data and let $\Tbar{X}$ be the tensor that is produced by the computed model. Then the
\emph{tensor completion score} (TCS) is
\begin{equation}
\label{eq:TCS}
  \text{TCS} = \frac{\| ({\bf 1} - \T{W}) \Hada (\T{X} - \Tbar{X}) \|}
  {\| ({\bf 1} - \T{W}) \Hada \T{X} \|}.
\end{equation}
In other words, the TCS is the relative error in the missing
entries. TCS is always nonnegative, and the best possible score is 0.

\subsection{Simulated data}
\label{sec:simulated}
We consider the performance of the methods on moderately-sized
problems of sizes $50 \times 40 \times 30$, $100 \times  80 \times
60$, and $150 \times 120 \times  90$. For all sizes, we set the number
of components in the CP model to be $R=5$. We test 60\%, 70\%, 80\%,
90\%, and 95\% missing data.
The experiments show that the underlying factors can be captured even
if the CP model
is fit to a tensor with a significant amount of missing data;
this is because the low-rank structure of the tensor
is being exploited. A rank-$R$ CP model for a tensor of size $I \times J \times K$
has $R(I+J+K-1)+1$ degrees of freedom.
The reason that the factors can be recovered accurately,
even with 95\% missing data, is that there is still a lot more data
than variables, i.e., the size of the data is equal to $0.05\, IJK$ which is
much greater than the $R(I+J+K-1)+1$ variables for large values of
$I$, $J$,and $K$. Because it
is a nonlinear problem,
we do not know exactly how many data entries are needed in order to
recover a CP model
of a low-rank tensor; however, a lower bound for the number of
entries needed in the matrix case has been derived in \cite{CaTa09}.

\subsubsection{Generating simulated data}
We create the test problems as follows.  Assume that the tensor size
is $I \times J \times K$ and that the number of factors is $R$.  We
generate factor matrices $\M{A}$, $\M{B}$, and $\M{C}$ of sizes $I
\times R$, $J \times R$, and $K \times R$, respectively, by randomly
choosing each entry from $\mathcal{N}(0,1)$ and then normalizing every
column to unit length.  Note that this method of generating the factor
matrices ensures that the solution is unique with probability one for
the sizes and number of components that we are considering because $R
\ll \min\{I,J,K\}$.\footnote{Since each matrix has $R$ linearly
independent columns with probability one and the $k$-rank (i.e., the
maximum value $k$ such that any 
$k$ columns are linearly independent) of each factor
matrix is $R$,  we satisfy the necessary conditions defined 
by Kruskal \cite{Kr77} for uniqueness.}

We then create the data tensor as
\begin{displaymath}
  \T{X} = \KOp{\M{A},\M{B},\M{C}} + \eta \;
  \frac{\norm{\T{X}}}{\norm{\T{N}}}\; \T{N}.
\end{displaymath}
Here $\T{N}$ is a noise tensor (of the same size as $\T{X}$) with
entries randomly selected from $\mathcal{N}(0,1)$, and $\eta$ is the
noise parameter. We use $\eta = 10\%$ in our experiments.

Finally, we set some entries of each generated tensor to be missing
according to an indicator tensor $\T{W}$.
We consider two situations for determining the missing data: randomly
missing entries and, in \App{appendix}, structured missing data in the
form of randomly
missing fibers. In the case of randomly missing entries, $\T{W}$ is a
binary tensor
such that exactly $\lfloor MIJK \rfloor$ randomly
selected entries are set to zero, where $M \in (0,1)$ defines the
percentage of missing data.
We require that every slice of $\T{W}$ (in every direction) have at
least one nonzero because otherwise we have no chance of recovering
the factor matrices;
this is related to the problem of \emph{coherence} in the matrix
completion problem where it is well-known that missing an entire row
(or a column)
of a matrix means that the true factors can never be recovered.
For each problem size and missing data percentage, thirty
independent test problems are generated.

\subsubsection{Results}
\label{sec:CI1}

In our experiments, CP-WOPT,
INDAFAC, and EM-ALS were indistinguishable in terms of
accuracy. \Fig{FMS_50-40-30} show box
plots\footnote{The box plots were generated using the \texttt{boxplot} command in the Statistics Toolbox of Matlab. The plots shown here use the ``compact'' plot style for that command. For more details on box plots and the use of this Matlab command, see \url{http://www.mathworks.com/access/helpdesk/help/toolbox/stats/boxplot.html}.}
of the FMS scores for problems of size $50 \times 40 \times 30$. Results for
problems of size  $100 \times
80 \times 60$ and $150 \times 120 \times 90$ are provided in
\App{random}. The plots summarize the FMS scores of the computed
factors for 30 independent problems for each problem size and
percentage of missing data: the median values are shown as black dots
inside a colored circle, the 25th and 75th percentile values
are indicated by the top and bottom edges of the solid bars extending from
the median values, and the outliers (corresponding to values more than
2.7 standard deviations from the means) are shown as
smaller colored circles. A lack of solid bars extending from the
median indicates that the 25th and 75th percentile values are equal
to the median.
The results shown are cumulative (i.e., the best so far) across the multiple starts, with
the first start using the $n$-mode singular vectors and the remaining starts using
random vectors.

For all methods, additional starting points improve performance; however,
no one method clearly requires fewer starting
points than the others. In general, using two or three
starting points suffices to solve all problems to high accuracy in
terms of FMS.
For problems with up to 90\% missing data, this is illustrated by 
median values which are all very close to one; further, the number of
outliers decreases with multiple starts. 

\FMSFIG{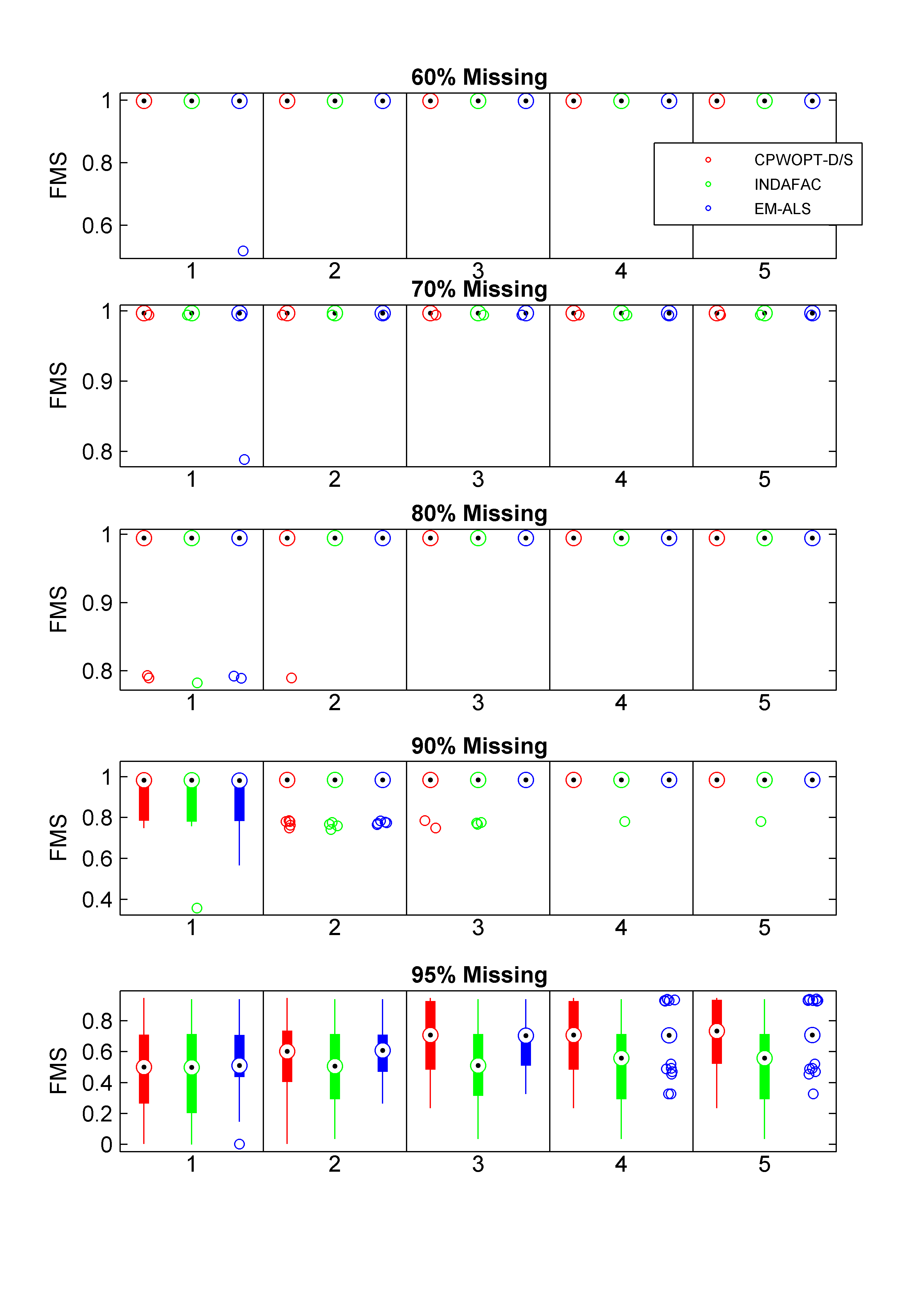}{50 \times 40 \times 30}{randomly missing entries}

Perhaps contrary to intuition, we note
that smaller problems are generally more difficult than larger
problems for the same percentage of missing data.
Thus, we can see in the bottom plot of \Fig{FMS_50-40-30} 
that the FMS scores remain low even with multiple
starts for problems of size $50 \times 40 \times 30$ with 95\%
missing data. It may be that some of these problems have gone below
the lower bound on the amount of data needed to find a solution; as
mentioned previously, a lower bound on how much data is needed is
known for the matrix case \cite{CaTa09} but finding such a bound is
still an open problem for higher-order tensors.
We can define the ratio $\rho$ as
\begin{equation}\label{eq:rho}
  \rho
  = \frac{\text{Number of known tensor entries}}
  {\text{Number of variables}}
  = \frac{(1-M)IJK}{R(I+J+K-1)+1},
\end{equation}
where smaller values of $\rho$ indicate more difficult problems.
For example, problems with 95\% missing data of size
$50 \times 40 \times 30$ ($\rho \approx 5 $) are more difficult than
those of size $100 \times 80 \times 60$ ($\rho \approx 20$), as illustrated
in \Fig{FMS_50-40-30} and \Fig{FMS_100-80-60}, respectively. In the
latter figure, corresponding to the larger size
(and thus larger value of $\rho$), nearly all problems are solved to high
accuracy by all methods.
Since this is a nonlinear problem, $\rho$ does not tell the
entire story, but it is at least a partial indicator of problem
difficulty.

The differentiator between the methods is computational time required. \Fig{Time} shows a
comparison of the sparse and dense versions of CPWOPT, along with INDAFAC and
EM-ALS. The timings reported in the figure are the sum of the times for all starting points.
The y-axis is time in seconds on a log scale. We present results for time rather than
iterations because the iterations for INDAFAC are much more expensive
than those for CPWOPT or EM-ALS.
We make
several observations.
For missing data levels up to
80\%,
the dense version of CPWOPT is
faster than INDAFAC, but EM-ALS is the fastest overall.
For 90\% and 95\% missing data, the sparse version of
CPWOPT is fastest (by more than a factor of 10 in some cases) with the exception
of the case of 90\% missing data for problems of size $150 \times 120
\times 90$. However, for this latter case, the difference is not
significant; furthermore, there are two problems where EM-ALS required
more than 10 times the median to find a correct factorization.

\TIMEFIG{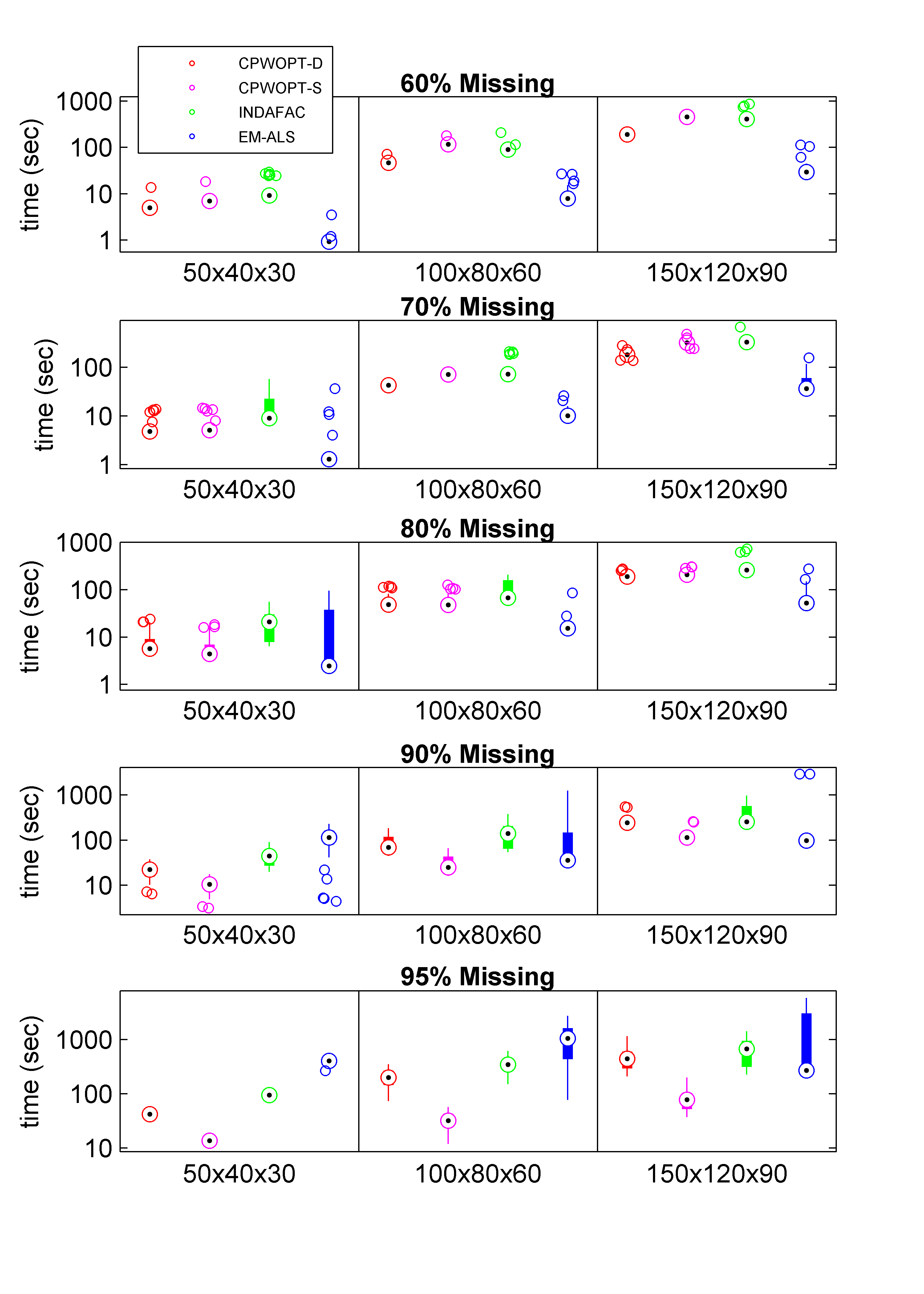}{randomly missing entries}

As demonstrated in earlier studies \cite{SDM10,ToBr05}, problems with
structured missing data are generally more difficult than those with randomly
missing values. We present results for problems of structured missing
data in \App{appendix}. The set-up is the same as we have used here,
except that we only go up to 90\% missing data because the accuracy of
all methods degrades significantly after that point. The
results indicate that accuracies of all methods start to diminish at
lower percentages of structured missing data than for those problems
with randomly missing data. Timings of the methods also indicate
similar behavior in the case of structure missing data: With 80\% or
less missing data EM-ALS is fastest, whereas the sparse version of
CP-WOPT is in general faster when there is 90\% missing data. 

Because it ignores missing values (random or structured), a major advantage of CP-WOPT is its
ability to perform sparse computations, enabling it to scale to very large
problems. Neither INDAFAC nor EM-ALS can do this.  The difficulty with
EM-ALS is that it must impute all missing values; therefore, it loses
any speed or scalability advantage of sparsity since it ultimately
operates on dense data. Although INDAFAC also ignores missing values,
its scalability is limited due to the expense of solving the Newton
equation. Even though the Newton system is extremely sparse, we can
see that the method is expensive even for moderate-sized problems such
as those discussed in this section. In the next section, we consider
very large problems, demonstrating this strength of CP-WOPT.

\subsection{Large-scale simulated data}
\label{sec:big}

A unique feature of CP-WOPT is that it can be applied to problems that
are too large to fit in memory if dense storage were
used. We consider two situations:
\begin{asparaenum}[(a)]
\item $500 \times 500 \times 500$ with 99\% missing data (1.25 million
  known values), and
\item $1000 \times 1000 \times 1000$ with 99.5\% missing data (5
                                million known values).
\end{asparaenum}

\subsubsection{Generating large-scale simulated data}

The method for generating test problems of this size is necessarily
different than the one for smaller problems presented in the previous section because we cannot, for
example, generate a full noise tensor. In fact, the tensors
$\T{Y} = \KOp{\M{A},\M{B},\M{C}}$ (the noise-free data) and $\T{N}$
(the noise) are never explicitly fully
formed for the large problems studied here.

Let the tensor size be $I \times J \times K$ and the missing value
rate be $M$.  We generate the factor matrices as described in
\Sec{simulated}, using $R=5$ components.  Next, we create the set
$\mathcal{I}$ with $(1-M)IJK$ randomly generated indices; this set
represents the indices of the known (non-missing) values in our tests.
The binary indicator tensor $\T{W}$ is stored as a sparse tensor
\cite{BaKo07} and defined by
\begin{displaymath}
  \TE{W}{ijk} =
  \begin{cases}
    1 & \text{if }(i,j,k) \in \mathcal{I}, \\
    0 & \text{otherwise}.
  \end{cases}
\end{displaymath}
Rather than explicitly forming all of $\T{Y} =
\KOp{\M{A},\M{B},\M{C}}$, we only calculate its values for those
indices in $\mathcal{I}$. This is analogous to the calculations
described for Step~2 of \Fig{sparse}. All missing entries of $\T{Y}$
are set to zero, and $\T{Y}$ is stored as a sparse tensor. Finally, we
set
\begin{displaymath}
  \T{X} = \T{Y} + \eta \;
  \frac{\norm{\T{X}}}{\norm{\T{N}}}\; \T{N},
\end{displaymath}
where $\T{N}$ is a \emph{sparse} noise tensor such that
\begin{displaymath}
  \TE{N}{ijk} =
  \begin{cases}
    \mathcal{N}(0,1) & \text{if }(i,j,k) \in \mathcal{I}, \\
    0 & \text{otherwise}.
  \end{cases}
\end{displaymath}
For each problem size, ten independent test problems are generated.
The initial guess is generated by computing the $n$-mode singular
vectors of the tensor with missing values filled in by zeros.

\subsubsection{Results for CP-WOPT on large-scale data}

\begin{figure}[htbp]
  \centering
  \subfloat[$500 \times 500 \times 500$ with $M=99\%$]%
  {\includegraphics[width=0.49\textwidth]{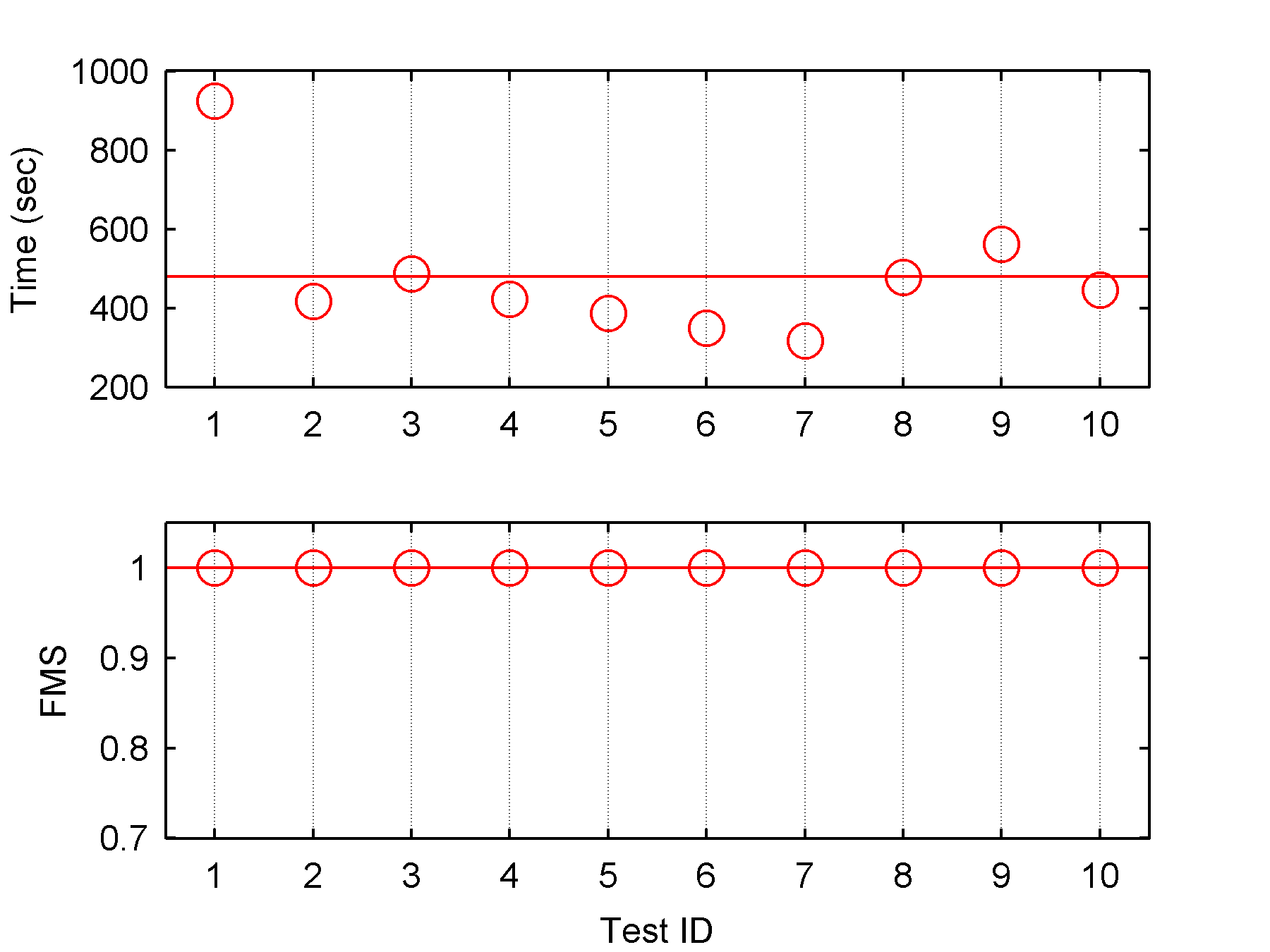}}
  \subfloat[$1000 \times 1000 \times 1000$ with $M=99.5\%$]%
  {\includegraphics[width=0.49\textwidth]{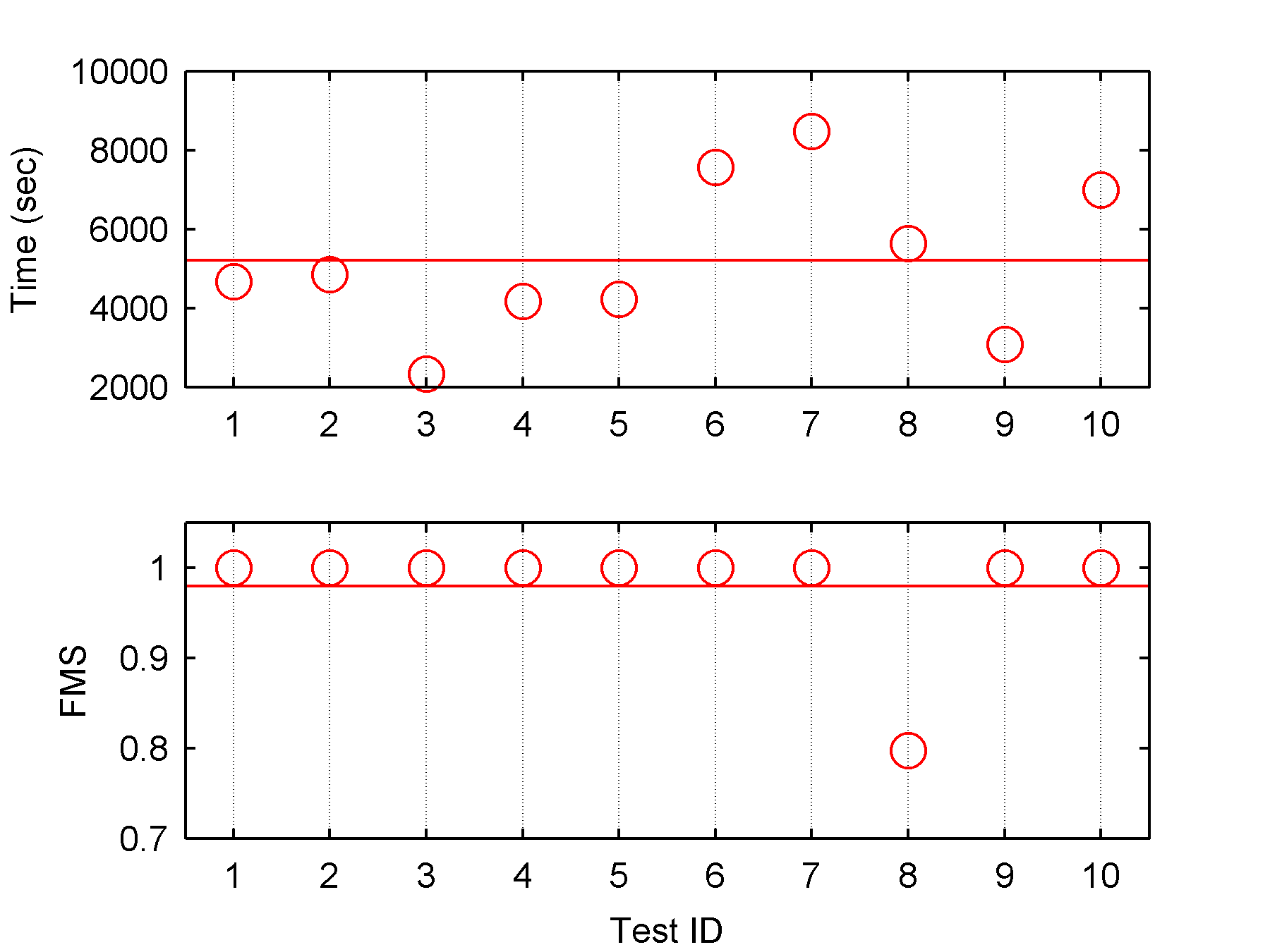}}
  \caption{Results of CP-WOPT-Sparse for large-scale problems with at
    least 99\% missing data. The means are shown as solid lines.}
  \label{fig:big}
\end{figure}

The computational set-up was the same as that described in \Sec{env}
except that the tolerance for the two-norm of the gradient divided by
the number of entries in the gradient was set to $10^{-10}$
(previously $10^{-8}$). The results across all ten runs for each
problem size are shown in \Fig{big}.

In the $500 \times 500 \times 500$ case, storing a dense version of
the tensor (assuming 8 bytes per entry) would require exactly 1GB.
Storing just 1\% of the data (1.25M entries) in sparse format
(assuming 32 bytes per
entry, 8 for the value and 8 for each of the indices) requires only
40MB.  In our randomly generated experiments, all ten problems were
solved with an FMS score greater than 0.99, with solve times ranging
between 200 and 1000 seconds.

In the $1000 \times 1000 \times 1000$ case, storing a dense version of
the tensor would require 8GB. Storing just 0.5\% of the data as a
sparse tensor (5M entries) requires 160MB. In our randomly generated
experiments, only nine of the ten problems were solved with an FMS
score greater than 0.99.
The solve times ranged from 2000 to 10000
seconds, approximately 10 times slower than the $500 \times 500 \times
500$ case, which had half as many variables and 1/4 the nonzero
entries.
In the failed case, the optimization method exited because the
relative change in the function tolerance was smaller than the
tolerance.
We restarted the optimization method with no changes except that we
use the solution that had
been computed previously as the initial guess. After 985 seconds of
additional work, an answer was found with an FMS score of 0.9999.


\subsection{EEG data}
\label{sec:EEG}
In this section, we demonstrate the use of CP-WOPT in multi-channel EEG (electroencephalogram) analysis
by capturing the underlying brain dynamics even in the presence of missing signals.
We use an EEG data set collected to observe the gamma activation during proprioceptive stimuli of left
and right hand \cite{MoHaAr07}. The data set contains
multi-channel signals (64 channels) recorded
from 14 subjects during stimulation of left and right hand (i.e., 28
measurements in total). For each measurement, the signal from each
channel is represented in both time and frequency domains using
a continuous wavelet transform and then vectorized (forming a vector of length 4392); in other words, each
measurement is represented by a \emph{channels} by
\emph{time-frequency} matrix. The data for all measurements are then
arranged as a \emph{channels} by \emph{time-frequency} by
\emph{measurements} tensor of size $64 \times 4392 \times 28$. For details about the data, see
\cite{MoHaAr07}.

\begin{figure}[htbp]
\centering
\subfloat[No missing entries]{\label{fig:EEGa}\includegraphics[width=.45\linewidth]{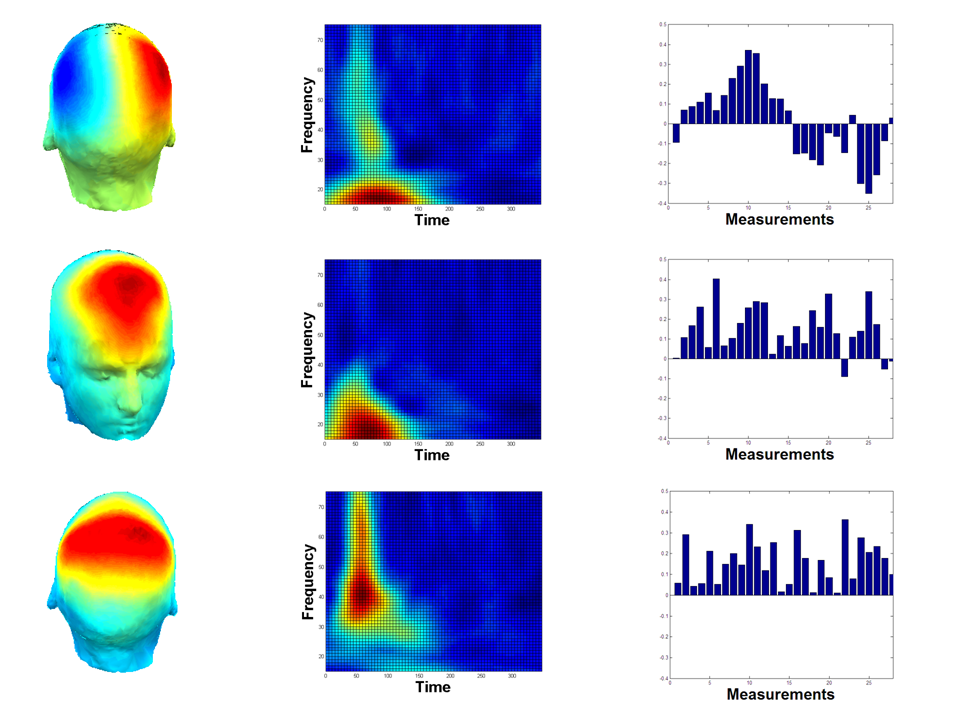}}
\subfloat[30 channels missing per measurement]{\label{fig:EEGb}\includegraphics[width=.45\linewidth]{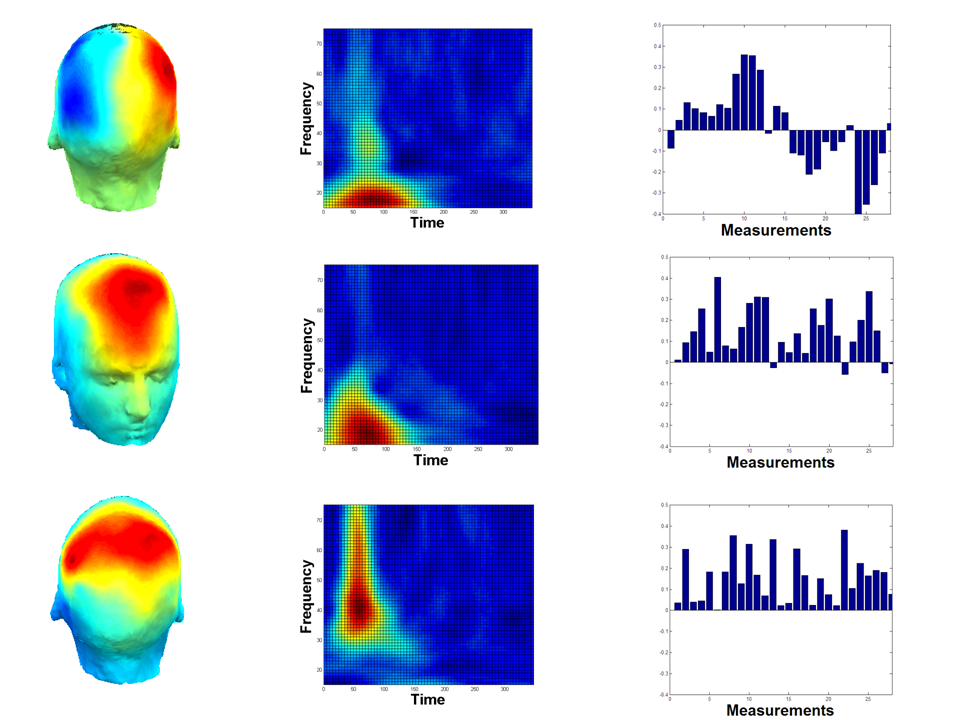}}
\caption[Illustration of CP factors for EEG data]{Columns of the CP factor matrices ($\M{A}$, $\M{B}$ and
  $\M{C}$ with $R=3$) extracted from the EEG data arranged as a
    \emph{channels} by \emph{time-frequency} by \emph{measurements}
    tensor.
    The 3-D head images were drawn using EEGLab~\cite{EEGLab}.} \label{fig:EEG}
\end{figure}

We model the data using a CP model with $R=3$, denoting
$\M{A}$, $\M{B}$, and $\M{C}$ as the extracted factor matrices corresponding to the channels,
time-frequency, and measurements modes, respectively.
We demonstrate the columns of the factor matrices in each mode in \Fig{EEGa}.
The 3-D head plots correspond to the columns of $\M{A}$, i.e., coefficients corresponding to the channels ranging
from low values in blue to high values in red. The time-frequency domain
representations correspond to the columns of $\M{B}$ rearranged as a matrix and again ranging from low values in blue to high values in red. The bar plots represent the columns of
$\M{C}$. Note that the three rows of images in \Fig{EEGa} (3-D head plot, matrix plot, bar plot) correspond to columns $r=1,2,3$ of the factor matrices ($\M{A}$, $\M{B}$, $\M{C}$), respectively. Observe that the first row of images highlights the differences between left and right hand stimulation
while the second and third rows of images pertain to frontal and parietal activations
that are shared by the stimuli. Unlike \cite{MoHaAr07}, we do not use
nonnegativity constraints; we convert tensor entries from complex
to real values by using the absolute values of the entries and center the data across the channels
mode before the analysis.

It is not uncommon in EEG analysis that the signals from some channels
are ignored due to malfunctioning of the electrodes. This corresponds to
missing fibers in the tensor (as in \Fig{fibers})
when we arrange the data as described above.
To reflect such cases of missing data, we randomly set data for one or more of the 64 channels for
each measurement to be missing, center the tensor across the channels mode ignoring the missing entries,
and then fit a CP model with $R=3$ to the resulting data using
CP-WOPT. Let $\Mbar{A}, \Mbar{B}, \Mbar{C}$ be the
factor matrices extracted from a tensor with missing entries using the CP-WOPT algorithm.

\Tab{EEG} presents the average FMS scores (over 50 instances of the
data with randomly missing channels) for different numbers of missing
channels per measurement.
As the number of
missing channels increases, the average FMS scores decrease as expected. However, even with up to 30
missing channels per measurement, or about $47\%$ of the data, the extracted factor matrices
match with the original factor matrices well,
with average FMS scores around $0.90$. \Fig{EEGb} presents the images for $\Mbar{A}, \Mbar{B}$
and $\Mbar{C}$ corresponding to those for $\M{A}, \M{B}, \M{C}$ in \Fig{EEGa}. We see that the
underlying brain dynamics are still captured even when $30$ channels per measurement are missing;
only slight local differences can be observed between the original data and the model generated using CP-WOPT.

It can be argued that the activations of the electrodes are highly correlated and
even if some of the electrodes are removed,
the underlying brain dynamics can still be captured. However, in these experiments
we do not set the same channels to missing for each measurement;
the channels are randomly missing from each measurement. There are cases when a CP
model may not be able to recover factors when data has certain
patterns of missing values, e.g., missing the signals from the same side of the brain
for all measurements. However, such data is not typical in practice.

\begin{table}[htpb]
  \caption[EEG Results]{\emph{EEG Analysis with Incomplete Data:} The
    similarity between the columns of $\M{A}, \M{B}, \M{C}$ and the
    columns of factor matrices extracted by CP-WOPT. The similarity is
    measured in terms of FMS defined in \Eqn{FMS}.}
  \label{tab:EEG}
  \begin{center}
  \begin{tabular}{|c | c |  }
    \hline
     Number of  &  CP-WOPT \\
     Missing Channels &  \\ \hline \hline
    $1$  &  $0.9959$ \\ \hline
    $10$ &  $0.9780$ \\ \hline
    $20$ &  $0.9478$ \\ \hline
    $30$ &  $0.8949$ \\ \hline
    $40$ &  $0.6459$ \\ \hline
  \end{tabular}
  \end{center}
\end{table}

\Tab{EEG} shows the FMS scores for varying amounts of missing
data. The solution does not seriously degrade until 40 of 64 channels
are removed.

\subsection{Network traffic data}

In the previous section, we were interested in tensor factorizations in the presence of missing data.
In this section, we address the problem of recovering missing entries of a tensor; in other words,
the \emph{tensor completion problem}. One application domain where this problem is frequently encountered is computer network traffic analysis.
Network traffic data consists of traffic matrices (TM), which record the amount of network data exchanged between source and
destination pairs (computers, routers, etc.). Since TMs evolve over time, network data can be represented as a tensor. For instance, G\'eant data \cite{UhQuLeBa06}
records the traffic exchanged between 23 routers over a period of 4 months collected using 15-minute intervals. This data set forms
a third-order tensor with \emph{source routers}, \emph{destination routers} and \emph{time} modes, and
each entry indicates the amount of traffic sent from a source to a destination during a particular time interval.
Missing data arises due to the expense of the data collection process.

\begin{figure*}[htbp]
\centering
{\includegraphics[width=3.5in]{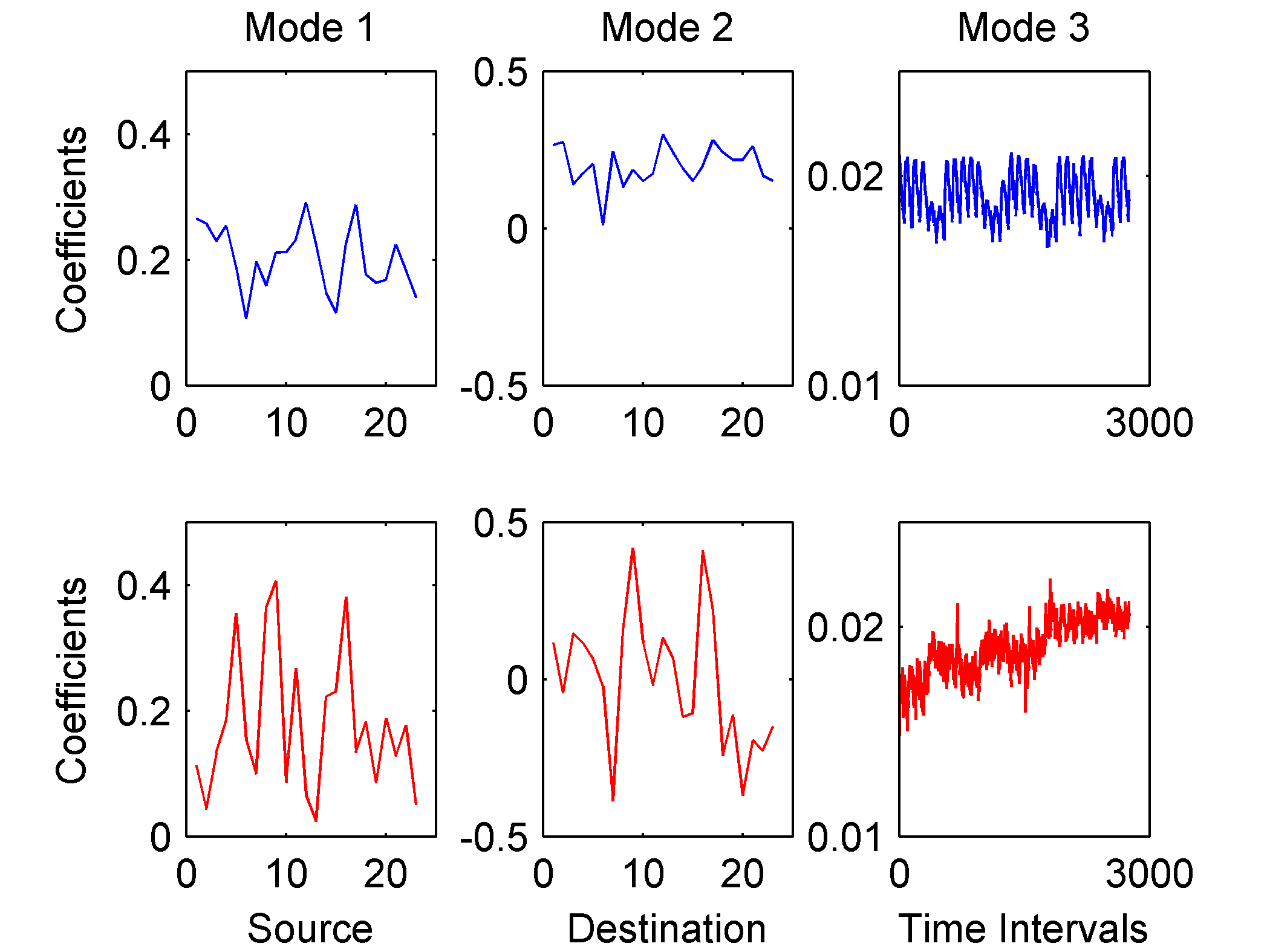}}
\caption[CP Factors]{Factor matrices extracted from G\'eant data using a 2-component CP model. The first row illustrates the first column of the factor matrices in each mode and the second row shows the second column of the factor matrices.} \label{fig:TMFactors}
\end{figure*}

In our study, we used the G\'eant data collected in April\footnote{The data was incomplete
for the full four-month period represented in the data, e.g., 6 days of data was missing at the end of February. We used a subset of the data corresponding to a period with no missing time slices.} stored as a tensor of size $23 \times 23 \times 2756$.
Let $\T{T}$ represent this raw data tensor. In order to adjust for scaling bias in the amount of traffic, we preprocess $\T{T}$ as $\T{X}=\log(\T{T}+1)$.
\Fig{TMFactors} presents the results of 2-component CP model (i.e., $R=2$);
the first and second row correspond the first and second column of the factor matrices of the model, respectively.
This CP model fits the data well but not perfectly and there is some unexplained variation left in the residuals.
Let $\That{X}$ be the tensor constructed using the computed factor matrices. The modeling error, defined as $\frac{\norm{\T{X}-\That{X}}}{\norm{\T{X}}}$, is approximately $0.31$ for the 2-component CP model computed. Even
though extracting more components slightly lowers the modeling error, models with more components do not look appropriate
for the data.

\begin{figure}[htbp]
\centering
{\includegraphics[width=2in]{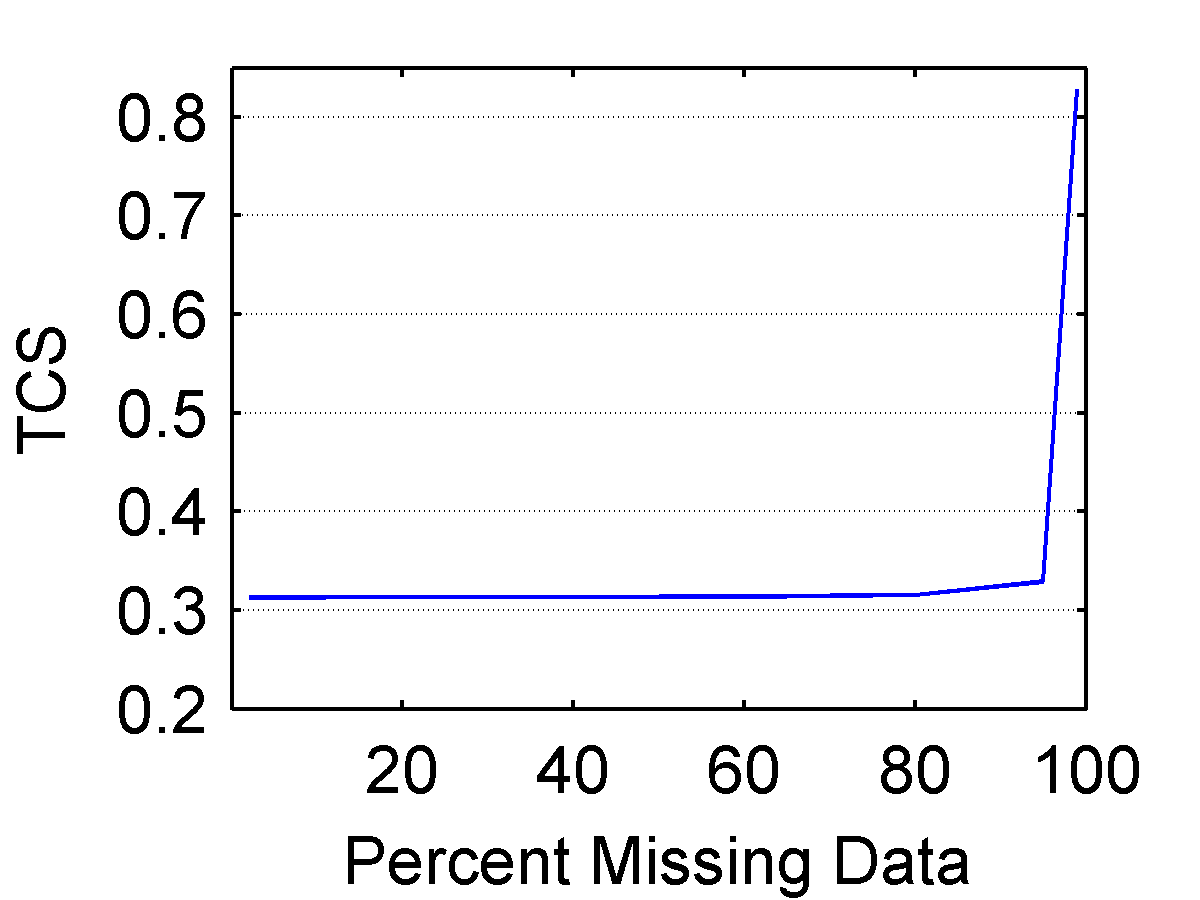}}
\caption[Tensor Completion Score]{Tensor Completion Score for different amounts of missing data for G\'eant data when a 2-component CP model fitted using CP-WOPT is used to recover the missing entries.} \label{fig:TCS}
\end{figure}

In order to assess the performance of CP-WOPT in terms of recovering missing data, randomly chosen entries of $\T{X}$
are set to missing and a 2-component CP model is fit to the altered data.
The extracted factor matrices are then used to reconstruct the data and fill in the missing entries. These recovered
values are compared with the actual values based on the tensor completion score (TCS) defined in \Eqn{TCS}.
\Fig{TCS} presents the average TCS (across 30 instances) for different amounts of missing data. We observe that the average TCS is around $0.31$
when there is little missing data and increases very slowly as we increase the amount of missing data. The average TCS
is only slightly higher, approximately $0.33$, even when $95\%$ of the entries are missing. However, the average TCS increases sharply when the amount of missing data is $99\%$. Note that even if there is no missing data, there will be completion error due to the modeling error, i.e., $0.31$. \Fig{TCS} demonstrates the robustness of CP-WOPT, illustrating that the TCS can be kept close to the level of the modeling error using this method, even when the amount of missing data is high.

We have addressed the tensor completion problem using a CP model, which gives easily-interpretable factors in each mode. However, for the tensor completion problem, the recovery of the missing entries is more important than the underlying factors and their interpretability. Therefore, modeling the data using a restricted model like CP may not be the best approach in this case. It may be possible to achieve lower modeling error using a more flexible model such as a Tucker model \cite{Tu66}. If this model is fit to the data using algorithms akin to CP-WOPT, missing entries may be recovered more accurately. This is a topic of future research.

\section{Conclusions}
\label{sec:conclusions}

The closely related problems of matrix factorizations with missing
data and matrix completion have recently been receiving a lot of
attention. In this paper, we consider the more general problem of
tensor factorization in the presence of missing data,
formulating the canonical
tensor decomposition for incomplete tensors as a weighted least
squares problem. Unlike imputation-based techniques,
this formulation ignores the missing entries and models only the known data entries.
We develop a scalable algorithm called CP-WOPT using gradient-based optimization
to solve the weighted least squares formulation of the CP problem. The
code has been released in version 2.4 of the Tensor Toolbox for MATLAB
\cite{TTB}.

Our numerical studies suggest that the proposed CP-WOPT approach
is accurate and scalable. CP-WOPT can recover the underlying factors
successfully with large amounts of missing data, e.g., 90\% missing entries for tensors
of size $50 \times 40 \times 30$. We have also studied how CP-WOPT can scale to problems of larger sizes, e.g.,
$1000 \times 1000 \times 1000$, and recover CP factors from large, sparse tensors with 99.5\% missing data.
Moreover, through the use of both dense and sparse implementations of the algorithm, we have shown
that CP-WOPT was always faster in our studies when compared to the best alternative approach
based on second-order optimization (INDAFAC) and even faster than
EM-ALS for high percentages of missing data.

We consider the practical use of CP-WOPT
algorithm in two different applications which demonstrate the
effectiveness of CP-WOPT even when there may be low correlation with a
multi-linear model. In multi-channel EEG analysis, the factors extracted by the
CP-WOPT algorithm can capture brain dynamics even if signals from some channels are missing, suggesting that
practitioners can now make better use of incomplete data in their
analyses. 
We note that the EEG data was centered by
using the means of only the known entries; however, robust techniques for
centering incomplete data are needed and is a topic of future investigation.
In network traffic analysis,
CP-WOPT algorithm can be used in the context of tensor completion and
recover the missing network traffic data. 

Although 90\% or more missing data may not seem practical, there are
situations where it is useful to purposely omit data.
For example, if the data is very large, excluding some data
will speed up the computation and enable it to fit in memory (for
large-scale data like the example of a $1000 \times 1000 \times 1000$
tensor that would normally require 8GB of storage). It is also common
to leave out data for the purpose of rank determination or model
assessment via cross-validation. For these reasons, it is useful to
consider high degrees
of missing data.

In future studies, we plan to extend our results in several
directions. We will include constraints such as non-negativity and penalties to encourage
sparsity, which enable us to find more meaningful latent factors from
large-scale sparse data. Finally, we will consider the problem of
collective factorizations with missing data, where we are jointly
factoring multiple tensors with shared factors.

\appendix
\section{Additional comparisons for randomly missing data}
\label{sec:random}

Figures \ref{fig:FMS_100-80-60} and \ref{fig:FMS_150-120-90} show box
plots of the FMS scores for problems of size $100 \times 80 \times 60$
and $150 \times 120 \times 90$, respectively. See \Sec{simulated} and \Fig{FMS_50-40-30} for a detailed description of the experimental setup and the results for analogous problems of size $50 \times 40 \times 30$. We can see from these
figures that the three methods being compared (CP-WOPT, INDAFAC,
EM-ALS) are indistinguishable in terms of accuracy, as was the case for the smaller sized problems. Furthermore, comparable improvements using multiple starts can be seen for these larger problems; we are able to solve a majority of the problems with very high accuracy using five or fewer starting points (the $n$-mode singular vectors used as the basis for the first run and random points used for the others).

\FMSFIG{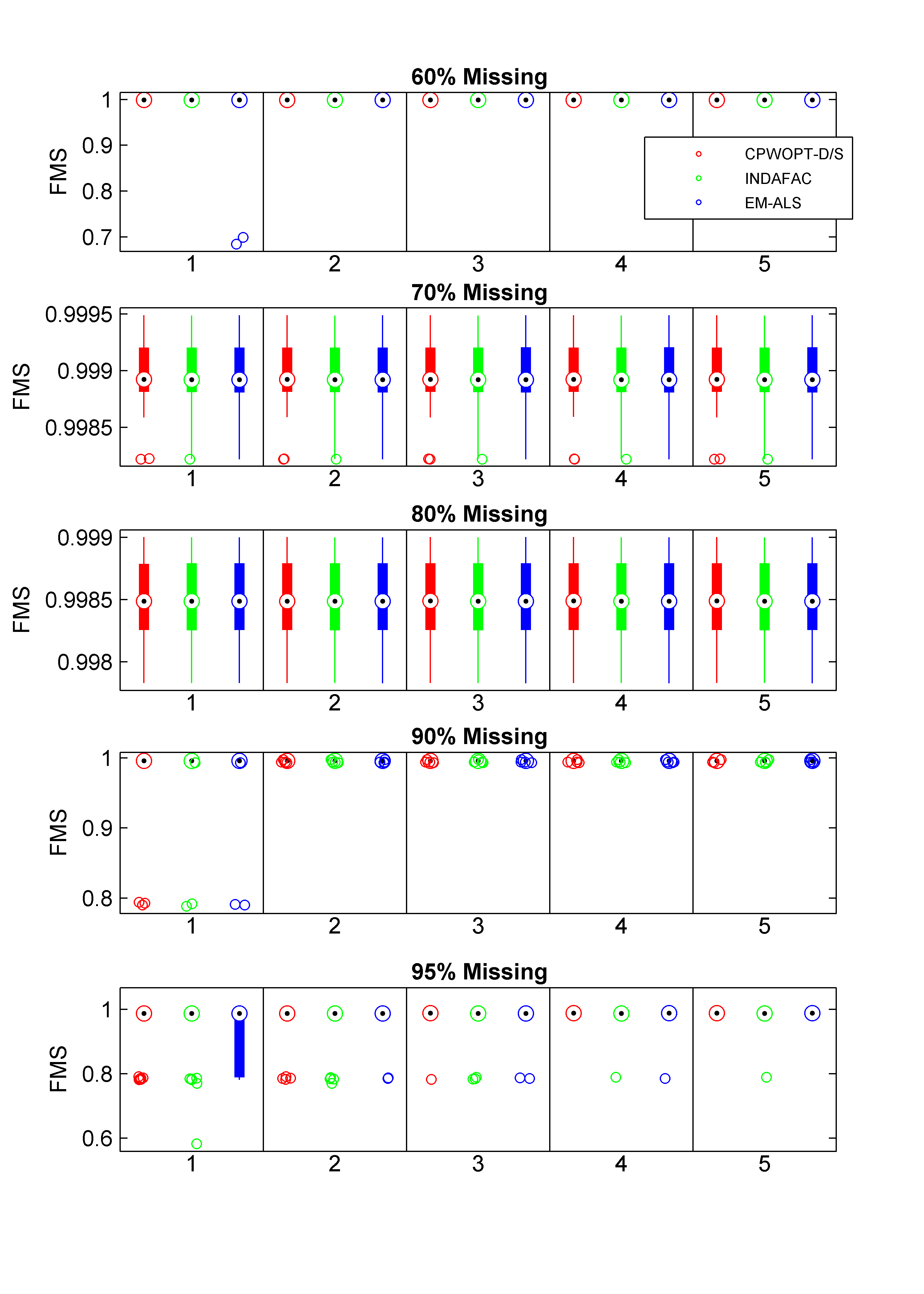}{100 \times 80 \times 60}{randomly missing entries}
\FMSFIG{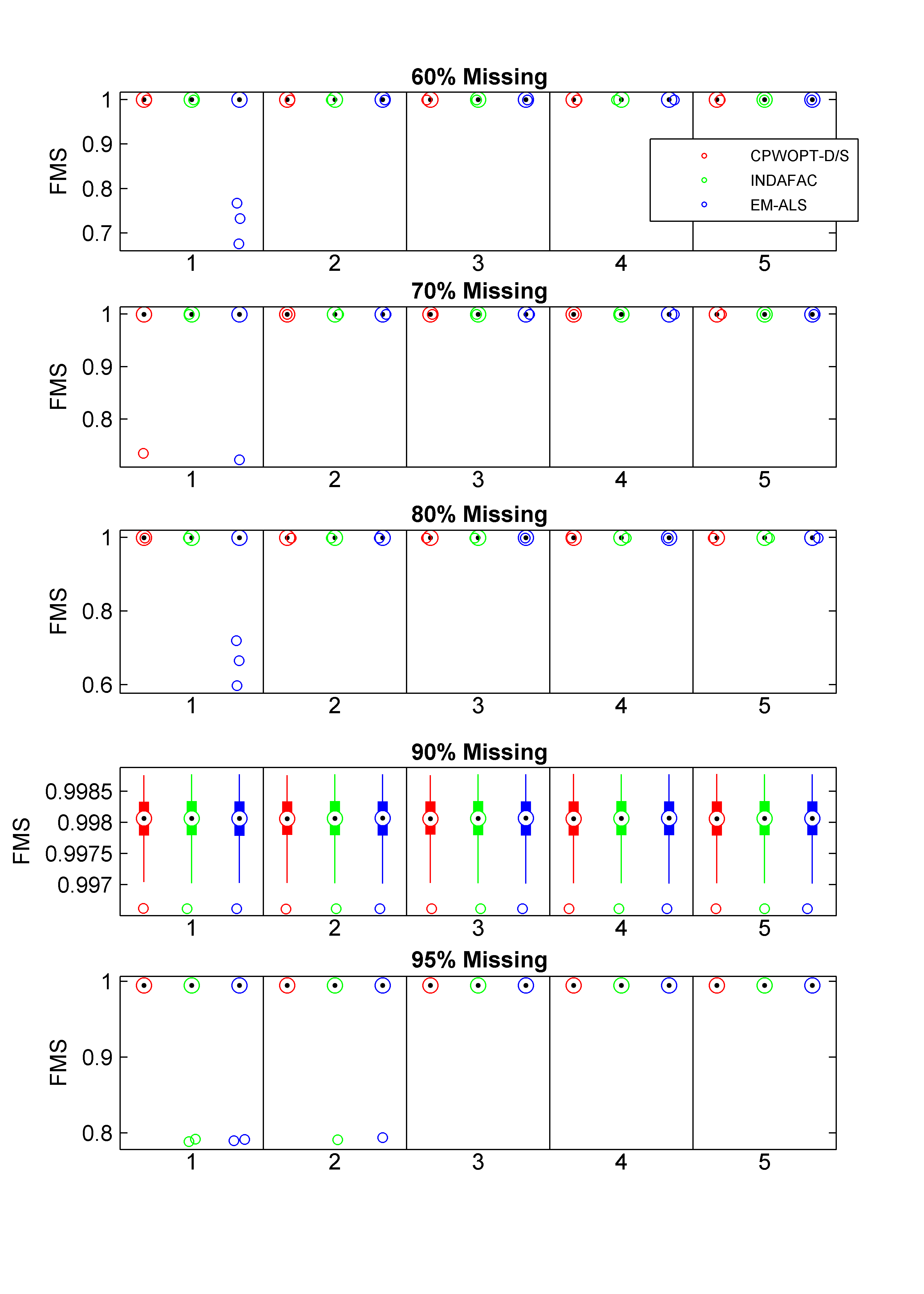}{150 \times 120 \times 90}{randomly missing entries}

\section{Comparisons for structured missing data}
\label{sec:appendix}

We use the same experimental set-up as described in \Sec{simulated}, 
with the exception that we performed experiments using problems with
up to $90\%$ missing data (as opposed to $95\%$) and we only show
results for three starting point (the first is the $n$-mode singular
vectors and the remaining two are random). 
In the case of randomly missing fibers, without loss of
generality, we consider missing fibers in the third mode only. In that
case, $\M{W}$ is an $I \times J$ binary matrix with exactly $\lfloor
MIJ \rfloor$ randomly selected entries are set to zero, and the binary
tensor $\T{W}$ is created by stacking $K$ copies of $\M{W}$ together.
We again require that every slice of $\T{W}$ (in every direction) have at
least one nonzero, which is equivalent to requiring that $\M{W}$ has
no zero rows or columns.
For each problem size and missing data percentage, thirty
independent test problems were generated.

The average FMS results are presented in Figures~\ref{fig:FMS_50-40-30_s}, \ref{fig:FMS_100-80-60_s} and \ref{fig:FMS_150-120-90_s} for problems of sizes $50 \times 40 \times 30$, $100 \times 80 \times 60$ and $150 \times 120 \times 90$, respectively. As for the problems with randomly missing data, we see that as the percentage of missing data increases, the average FMS scores tend to decrease. However, one notable difference is that the average FMS scores for structured missing data are in general lower than those for problems with comparable amounts of randomly missing data. This indicates that problems with structured missing data are more difficult than those with randomly missing data.

\FMSFIG{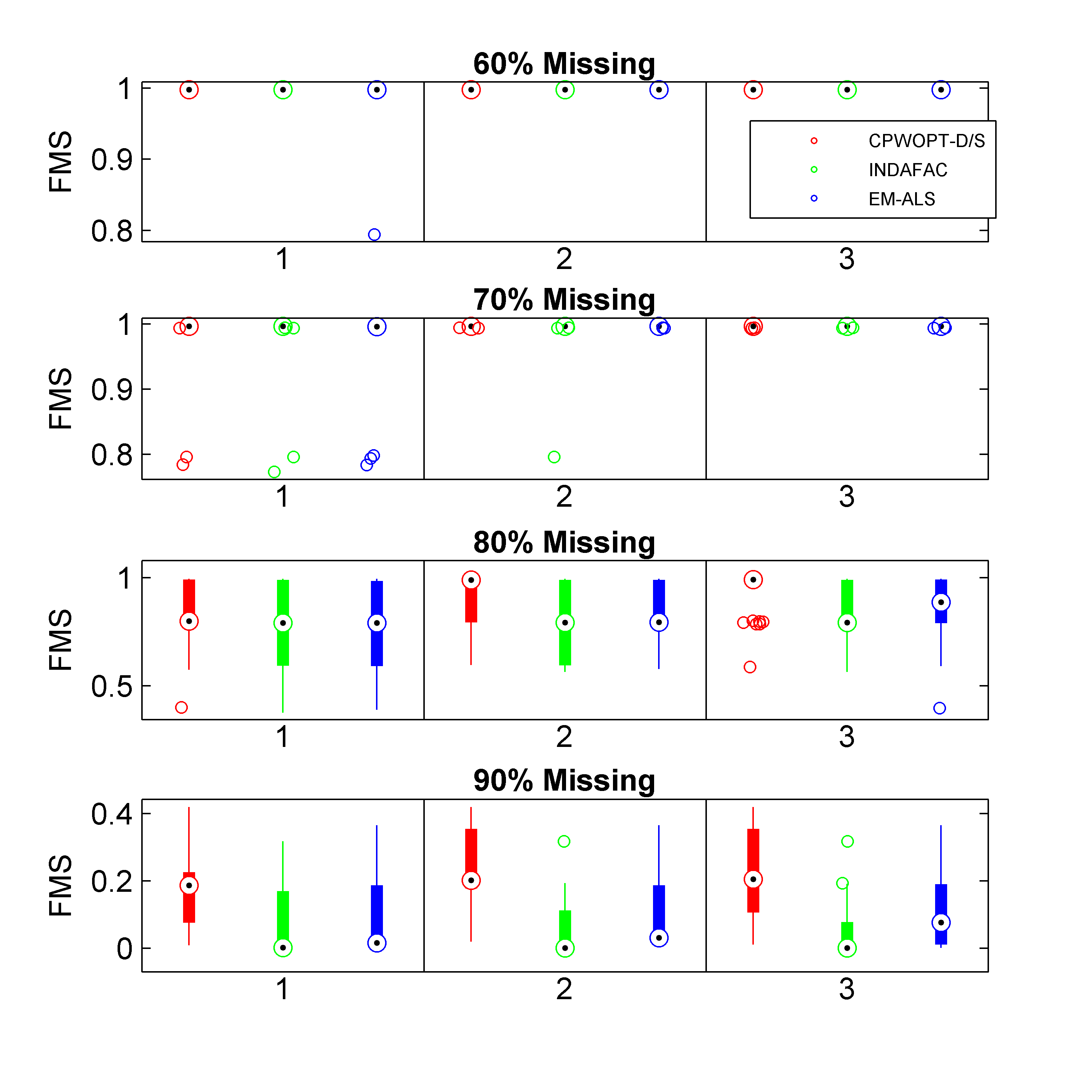}{50 \times 40 \times 30}{randomly missing fibers}
\FMSFIG{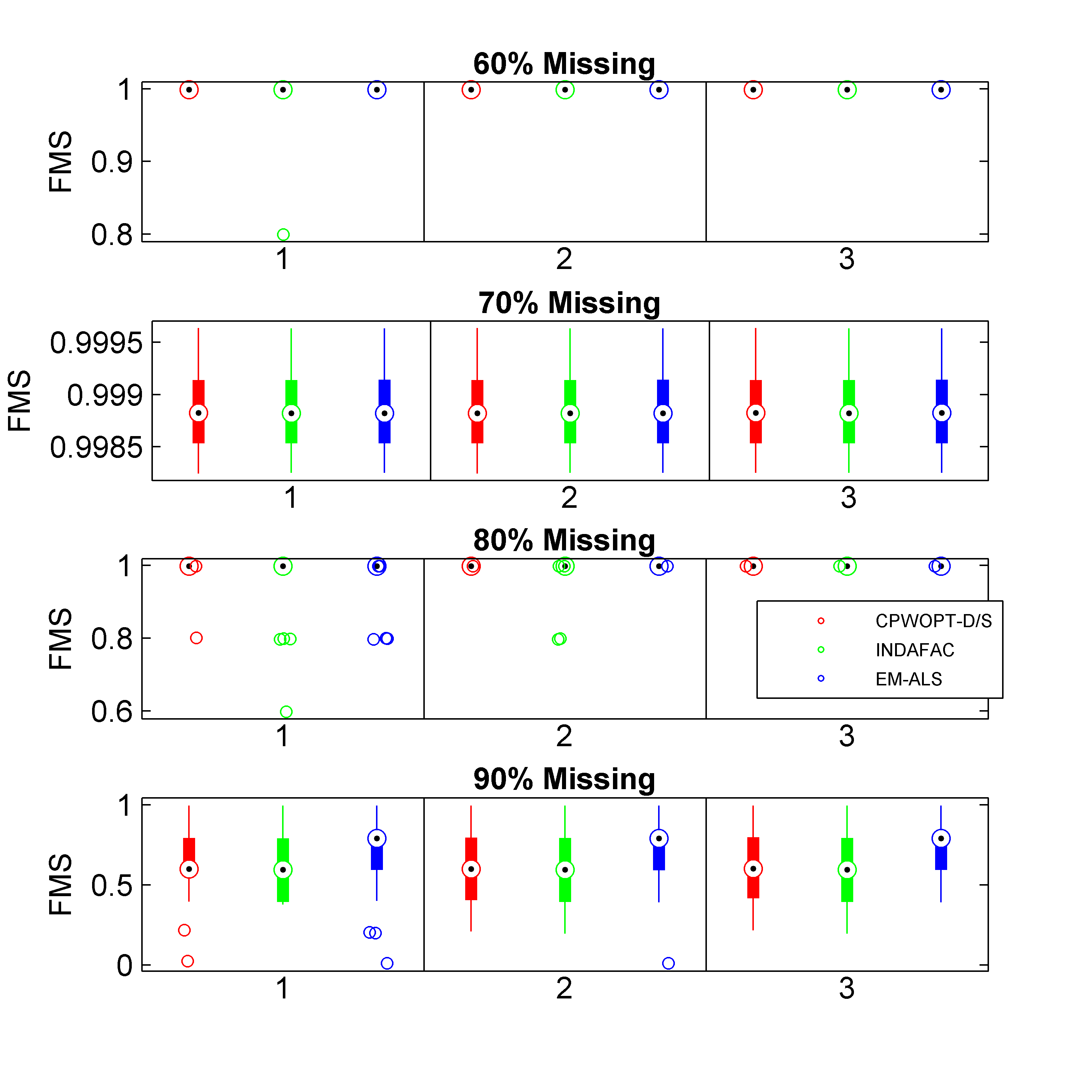}{100 \times 80 \times 60}{randomly missing fibers}
\FMSFIG{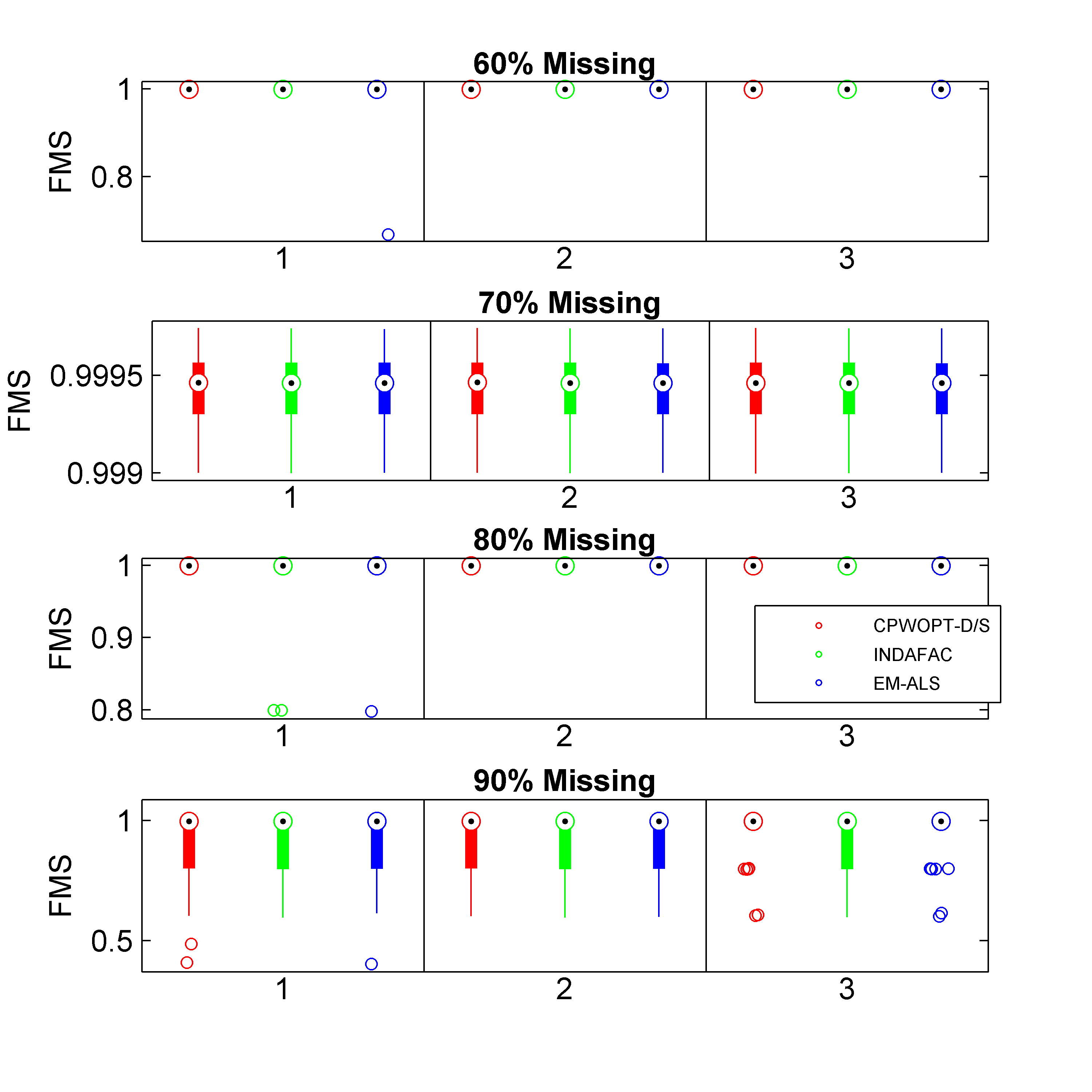}{150 \times 120 \times 90}{randomly missing fibers}

The computational times required to compute the CP models using the different methods are present in \Fig{Time_s}, where we observe comparable results to those presented in \Fig{Time} for the problems with randomly missing data. With 80\% or less structured missing data EM-ALS is 
fastest, whereas the sparse version of CP-WOPT is in general faster when there is more than 80\% missing data.

\TIMEFIG{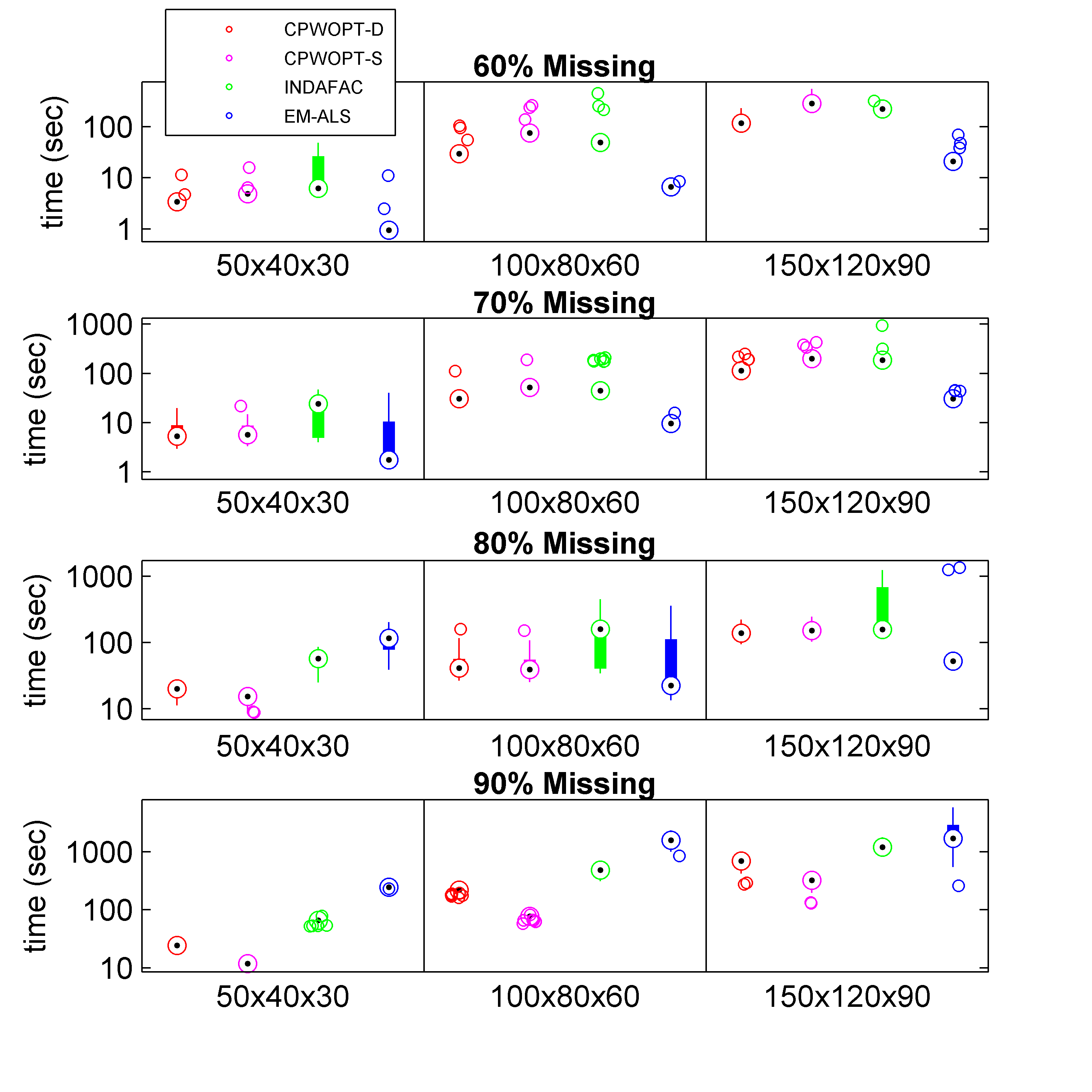}{randomly missing fibers}

\section*{Acknowledgments}
We thank David Gleich and the anonymous referees for helpful comments
which greatly improved the presentation of this manuscript.




\begin{thebibliography}{10}
\expandafter\ifx\csname url\endcsname\relax
  \def\url#1{\texttt{#1}}\fi
\expandafter\ifx\csname urlprefix\endcsname\relax\def\urlprefix{URL }\fi
\expandafter\ifx\csname href\endcsname\relax
  \def\href#1#2{#2} \def\path#1{#1}\fi

\bibitem{SDM10}
E.~Acar, D.~M. Dunlavy, T.~G. Kolda, M.~M{\o}rup,
  \href{http://www.siam.org/proceedings/datamining/2010/dm10_061_acare.pdf}{Sc%
alable tensor factorizations with missing data}, in: Proceedings of the Tenth
  SIAM International Conference on Data Mining, SIAM, 2010, pp. 701--712.
\newline\urlprefix\url{http://www.siam.org/proceedings/datamining/2010/dm10_06%
1_acare.pdf}

\bibitem{ZhRoWiQi09}
Y.~Zhang, M.~Roughan, W.~Willinger, L.~Qiu, Spatio-temporal compressive sensing
  and internet traffic matrices, in: SIGCOMM '09: Proceedings of the ACM
  SIGCOMM 2009 conference on Data communication, ACM, New York, NY, USA, 2009,
  pp. 267--278.
\newblock \href {http://dx.doi.org/10.1145/1592568.1592600}
  {\path{doi:10.1145/1592568.1592600}}.

\bibitem{BuFi05}
A.~M. Buchanan, A.~W. Fitzgibbon, Damped {N}ewton algorithms for matrix
  factorization with missing data, in: CVPR'05: 2005 IEEE Computer Society
  Conference on Computer Vision and Pattern Recognition, Vol.~2, IEEE Computer
  Society, 2005, pp. 316--322.
\newblock \href {http://dx.doi.org/10.1109/CVPR.2005.118}
  {\path{doi:10.1109/CVPR.2005.118}}.

\bibitem{AcYe09}
E.~Acar, B.~Yener, Unsupervised multiway data analysis: A literature survey,
  IEEE Transactions on Knowledge and Data Engineering 21~(1) (2009) 6--20.
\newblock \href {http://dx.doi.org/10.1109/TKDE.2008.112}
  {\path{doi:10.1109/TKDE.2008.112}}.

\bibitem{KoBa09}
T.~G. Kolda, B.~W. Bader, Tensor decompositions and applications, SIAM Review
  51~(3) (2009) 455--500.
\newblock \href {http://dx.doi.org/10.1137/07070111X}
  {\path{doi:10.1137/07070111X}}.

\bibitem{MiMaVaNi04}
F.~Miwakeichi, E.~Mart\'{i}nez-Montes, P.~A. Valdés-Sosa, N.~Nishiyama,
  H.~Mizuhara, Y.~Yamaguchi, Decomposing {EEG} data into space-time-frequency
  components using parallel factor analysis, NeuroImage 22~(3) (2004)
  1035--1045.
\newblock \href {http://dx.doi.org/10.1016/j.neuroimage.2004.03.039}
  {\path{doi:10.1016/j.neuroimage.2004.03.039}}.

\bibitem{ToBr05}
G.~Tomasi, R.~Bro, {PARAFAC} and missing values, Chemometrics and Intelligent
  Laboratory Systems 75~(2) (2005) 163--180.
\newblock \href {http://dx.doi.org/10.1016/j.chemolab.2004.07.003}
  {\path{doi:10.1016/j.chemolab.2004.07.003}}.

\bibitem{Orekhov2003}
V.~Y. Orekhov, I.~Ibraghimov, M.~Billeter, Optimizing resolution in
  multidimensional {NMR} by three-way decomposition, Journal of Biomolecular
  NMR 27 (2003) 165--173.
\newblock \href {http://dx.doi.org/10.1023/A:1024944720653}
  {\path{doi:10.1023/A:1024944720653}}.

\bibitem{Geng2009}
X.~Geng, K.~Smith-Miles, Z.-H. Zhou, L.~Wang, Face image modeling by
  multilinear subspace analysis with missing values, in: MM '09: Proceedings of
  the seventeen ACM international conference on Multimedia, ACM, 2009, pp.
  629--632.
\newblock \href {http://dx.doi.org/10.1145/1631272.1631373}
  {\path{doi:10.1145/1631272.1631373}}.

\bibitem{Ru74}
A.~Ruhe, Numerical computation of principal components when several
  observations are missing, Tech. Rep. UMINF-48-74, Department of Information
  Processing, Institute of Mathematics and Statistics, University of Umea,
  Umea, Sweden (1974).

\bibitem{GaZa79}
K.~R. Gabriel, S.~Zamir, \href{http://www.jstor.org/stable/1268288}{Lower rank
  approximation of matrices by least squares approximation with any choice of
  weights}, Technometrics 21~(4) (1979) 489--498.
\newline\urlprefix\url{http://www.jstor.org/stable/1268288}

\bibitem{CaTa09}
E.~J. Cand{\`e}s, T.~Tao, \href{http://arxiv.org/abs/0903.1476}{The power of
  convex relaxation: Near-optimal matrix completion}, arXiv:0903.1476v1 (Mar.
  2009).
\newline\urlprefix\url{http://arxiv.org/abs/0903.1476}

\bibitem{CaPl09}
E.~J. Cand{\`e}s, Y.~Plan, \href{http://arxiv.org/abs/0903.3131}{Matrix
  completion with noise}, arXiv:0903.3131v1 (Mar. 2009).
\newline\urlprefix\url{http://arxiv.org/abs/0903.3131}

\bibitem{CaCh70}
J.~D. Carroll, J.~J. Chang, Analysis of individual differences in
  multidimensional scaling via an {N}-way generalization of ``{Eckart-Young}''
  decomposition, Psychometrika 35 (1970) 283--319.
\newblock \href {http://dx.doi.org/10.1007/BF02310791}
  {\path{doi:10.1007/BF02310791}}.

\bibitem{Ha70}
R.~A. Harshman, Foundations of the {PARAFAC} procedure: Models and conditions
  for an ``explanatory" multi-modal factor analysis, UCLA working papers in
  phonetics 16 (1970) 1--84, available at
  \url{http://www.psychology.uwo.ca/faculty/harshman/wpppfac0.pdf}.

\bibitem{KoBaKe05}
T.~G. Kolda, B.~W. Bader, J.~P. Kenny, Higher-order web link analysis using
  multilinear algebra, in: ICDM 2005: Proceedings of the 5th IEEE International
  Conference on Data Mining, IEEE Computer Society, 2005, pp. 242--249.
\newblock \href {http://dx.doi.org/10.1109/ICDM.2005.77}
  {\path{doi:10.1109/ICDM.2005.77}}.

\bibitem{Br06a}
R.~Bro, Review on multiway analysis in chemistry---2000--2005, Critical Reviews
  in Analytical Chemistry 36~(3--4) (2006) 279--293.
\newblock \href {http://dx.doi.org/10.1080/10408340600969965}
  {\path{doi:10.1080/10408340600969965}}.

\bibitem{AcBiBiBr07}
E.~Acar, C.~A. Bingol, H.~Bingol, R.~Bro, B.~Yener, Multiway analysis of
  epilepsy tensors, Bioinformatics 23~(13) (2007) i10--i18.
\newblock \href {http://dx.doi.org/10.1093/bioinformatics/btm210}
  {\path{doi:10.1093/bioinformatics/btm210}}.

\bibitem{MoHaAr07}
M.~M{\o}rup, L.~K. Hansen, S.~M. Arnfred, {ERPWAVELAB} a toolbox for
  multi-channel analysis of time-frequency transformed event related
  potentials, Journal of Neuroscience Methods 161~(2) (2007) 361--368.
\newblock \href {http://dx.doi.org/10.1016/j.jneumeth.2006.11.008}
  {\path{doi:10.1016/j.jneumeth.2006.11.008}}.

\bibitem{AcKoDu09}
E.~Acar, T.~Kolda, D.~Dunlavy, An optimization approach for fitting canonical
  tensor decompositions, Tech. Rep. SAND2009-0857, Sandia National
  Laboratories, Albuquerque, New Mexico and Livermore, California (2009).

\bibitem{DeLaRu77}
A.~P. Dempster, N.~M. Laird, D.~B. Rubin, Maximum likelihood from incomplete
  data via the em algorithm, Journal of the Royal Statistical Society, Series B
  39~(1) (1977) 1--38.

\bibitem{SrJa03}
N.~Srebro, T.~Jaakkola, Weighted low-rank approximations, in: IMCL-2003:
  Proceedings of the Twentieth International Conference on Machine Learning,
  2003, pp. 720--727.

\bibitem{Ki97}
H.~A.~L. Kiers, Weighted least squares fitting using ordinary least squares
  algorithms, Psychometrika 62~(2) (1997) 215--266.
\newblock \href {http://dx.doi.org/10.1007/BF02295279}
  {\path{doi:10.1007/BF02295279}}.

\bibitem{Br98}
R.~Bro, Multi-way analysis in the food industry: Models, algorithms, and
  applications, Ph.D. thesis, University of Amsterdam, available at
  \url{http://www.models.kvl.dk/research/theses/} (1998).

\bibitem{WaMa01}
B.~Walczak, D.~L. Massart, Dealing with missing data: Part {I}, Chemometrics
  and Intelligent Laboratory Systems 58~(1) (2001) 15--27.
\newblock \href {http://dx.doi.org/10.1016/S0169-7439(01)00131-9}
  {\path{doi:10.1016/S0169-7439(01)00131-9}}.

\bibitem{Pa97}
P.~Paatero, A weighted non-negative least squares algorithm for three-way
  ``{PARAFAC}'' factor analysis, Chemometrics and Intelligent Laboratory
  Systems 38~(2) (1997) 223--242.
\newblock \href {http://dx.doi.org/10.1016/S0169-7439(97)00031-2}
  {\path{doi:10.1016/S0169-7439(97)00031-2}}.

\bibitem{NoWr99}
J.~Nocedal, S.~J. Wright, Numerical Optimization, Springer, 1999.

\bibitem{BaKo07}
B.~W. Bader, T.~G. Kolda, Efficient {MATLAB} computations with sparse and
  factored tensors, SIAM Journal on Scientific Computing 30~(1) (2007)
  205--231.
\newblock \href {http://dx.doi.org/10.1137/060676489}
  {\path{doi:10.1137/060676489}}.

\bibitem{TTB}
B.~W. Bader, T.~G. Kolda, {MATLAB} tensor toolbox version 2.4,
  \url{http://csmr.ca.sandia.gov/~tgkolda/TensorToolbox/} (last accessed March,
  2010).

\bibitem{MoTh94}
J.~J. Mor\'{e}, D.~J. Thuente, Line search algorithms with guaranteed
  sufficient decrease, ACM Transactions on Mathematical Software 20~(3) (1994)
  286--307.
\newblock \href {http://dx.doi.org/10.1145/192115.192132}
  {\path{doi:10.1145/192115.192132}}.

\bibitem{DuKoAc10}
D.~M. Dunlavy, T.~G. Kolda, E.~Acar, Poblano v1.0: A {Matlab} toolbox for
  gradient-based optimization, Tech. Rep. SAND2010-1422, Sandia National
  Laboratories, Albuquerque, NM and Livermore, CA (Mar. 2010).

\bibitem{AnBr00}
C.~A. Andersson, R.~Bro, The {N}-way toolbox for {MATLAB}, Chemometrics and
  Intelligent Laboratory Systems 52~(1) (2000) 1--4, see also
  \url{http://www.models.kvl.dk/source/nwaytoolbox/}.
\newblock \href {http://dx.doi.org/10.1016/S0169-7439(00)00071-X}
  {\path{doi:10.1016/S0169-7439(00)00071-X}}.

\bibitem{INDAFAC}
G.~Tomasi, Incomplete data {PARAFAC} ({INDAFAC}),
  \url{http://www.models.kvl.dk/source/indafac/index.asp} (last accessed May,
  2009).

\bibitem{Kr77}
J.~B. Kruskal, Three-way arrays: rank and uniqueness of trilinear
  decompositions, with application to arithmetic complexity and statistics,
  Linear Algebra and its Applications 18~(2) (1977) 95--138.
\newblock \href {http://dx.doi.org/10.1016/0024-3795(77)90069-6}
  {\path{doi:10.1016/0024-3795(77)90069-6}}.

\bibitem{EEGLab}
A.~Delorme, S.~Makeig, {EEGLAB}: An open source toolbox for analysis of
  single-trial {EEG} dynamics, J. Neurosci. Meth. 134 (2004) 9--21.
\newblock \href {http://dx.doi.org/10.1016/j.jneumeth.2003.10.009}
  {\path{doi:10.1016/j.jneumeth.2003.10.009}}.

\bibitem{UhQuLeBa06}
S.~Uhlig, B.~Quoitin, J.~Lepropre, S.~Balon, Providing public intradomain
  traffic matrices to the research community, ACM SIGCOMM Computer
  Communication Review 36~(1) (2006) 83--86.
\newblock \href {http://dx.doi.org/10.1145/1111322.1111341}
  {\path{doi:10.1145/1111322.1111341}}.

\bibitem{Tu66}
L.~R. Tucker, Some mathematical notes on three-mode factor analysis,
  Psychometrika 31 (1966) 279--311.
\newblock \href {http://dx.doi.org/10.1007/BF02289464}
  {\path{doi:10.1007/BF02289464}}.

\end{thebibliography}
\end{document}